\documentclass[11pt]{amsart}
\usepackage{amssymb,amscd}

\newdimen\theight
\def\TeXref#1{%
              \leavevmode\vadjust{\setbox0=\hbox{{\tt
                      \quad\quad  {\small \textrm #1}}}%
              \theight=\ht0
              \advance\theight by \lineskip
              \kern -\theight \vbox to
              \theight{\rightline{\rlap{\box0}}%
              \vss}%
              }}%

\theoremstyle{plain}

\newtheorem{thm}{Theorem}[section]
\newtheorem{lem}[thm]{Lemma}
\newtheorem{cor}[thm]{Corollary}
\newtheorem{prop}[thm]{Proposition}

\theoremstyle{definition}

\newtheorem{conj}[thm]{Conjecture}

\theoremstyle{remark}

\newtheorem{claim}{Claim}

\newcommand{\AAA}{\text{$\mathcal{A}$}}

\newcommand{\DD}{\text{$\mathcal{D}$}}

\newcommand{\FF}{\mathcal{F}}
\newcommand{\GG}{\text{$\mathcal{G}$}}
\newcommand{\HH}{\text{$\mathcal{H}$}}

\newcommand{\OO}{\text{$\mathcal{O}$}}

\newcommand{\SSS}{\text{$\mathcal{S}$}}

\newcommand{\VV}{\text{$\mathcal{V}$}}
\newcommand{\WW}{\text{$\mathcal{W}$}}

          \newcommand{\bfD}{{\mathbf D}}

\newcommand{\bfHH}{\boldsymbol{\mathcal{H}}}

\newcommand{\bfDelta}{\boldsymbol{\Delta}}

\newcommand{\bfPi}{\boldsymbol{\Pi}}

\newcommand{\bfOmega}{\boldsymbol{\Omega}}

        \newcommand{\field}[1]{\text{$\mathbb{#1}$}}
        
        \newcommand{\Z}{\field{Z}}
        
        \newcommand{\R}{\field{R}}
        \newcommand{\C}{\field{C}}

\newcommand{\supp}{\operatorname{supp}}
\newcommand{\codim}{\operatorname{codim}}

\newcommand{\id}{\operatorname{id}}
\newcommand{\Tr}{\operatorname{Tr}}
\newcommand{\sign}{\operatorname{sign}}
\newcommand{\Fix}{\operatorname{Fix}}
\newcommand{\vol}{\operatorname{vol}}
\newcommand{\ev}{\operatorname{ev}}
\newcommand{\Aut}{\operatorname{Aut}}
\newcommand{\pr}{\operatorname{pr}}
\newcommand{\Diff}{\operatorname{Diff}}
\newcommand{\Ad}{\operatorname{Ad}}
\newcommand{\Hom}{\operatorname{Hom}}
\newcommand{\End}{\operatorname{End}}
\newcommand{\Pf}{\operatorname{Pf}}
\newcommand{\Pen}{\operatorname{Pen}}

\newcommand{\cinf}{\text{$C^\infty$}}
\newcommand{\cinfc}{\text{$C^\infty_c$}}

\newcommand{\D}{\text{$\Delta$}}

\begin{document}

\bibliographystyle{plain}


\title{Lefschetz distribution of Lie foliations}
\author[J.A. \'Alvarez L\'opez]{Jes\'us A. \'Alvarez L\'opez}
\address{Departamento de Xeometr\'{\i}a e Topolox\'{\i}a\\
         Facultade de Matem\'aticas\\
         Universidade de Santiago de Compostela\\
         15782 Santiago de Compostela\\ Spain}
\email{jalvarez@usc.es}
\thanks{Partially supported by MEC (Spain), grant MTM2004-08214}

\author[Y.A. Kordyukov]{Yuri A. Kordyukov}
\address{Institute of Mathematics\\
         Russian Academy of Sciences\\
         112~Chernyshevsky str.\\ 450077 Ufa\\ Russia}
\email{yurikor@matem.anrb.ru}
\thanks{Partially supported by the Russian Foundation of Basic Research
(grant no. 06-01-00208)}



\maketitle


\begin{abstract}
Let $\FF$ be a Lie foliation on a closed manifold $M$ with
structural Lie group $G$. Its transverse Lie structure can be
considered as a transverse action $\Phi$ of $G$ on $(M,\FF)$; {\it
i.e.\/}, an ``action'' which is defined up to leafwise homotopies.
This $\Phi$ induces an action $\Phi^*$ of $G$ on the reduced
leafwise cohomology $\overline H(\FF)$. By using leafwise Hodge
theory, the supertrace of $\Phi^*$ can be defined as a distribution
$L_{\text{\rm dis}}(\FF)$ on $G$ called the Lefschetz distribution
of $\FF$. A distributional version of the Gauss-Bonett theorem is
proved, which describes $L_{\text{\rm dis}}(\FF)$ around the
identity element. On any small enough open subset of $G$,
$L_{\text{\rm dis}}(\FF)$ is described by a distributional version
of the Lefschetz trace formula.
\end{abstract}

\tableofcontents

\section{Introduction}

Let $\FF$ be a \cinf\ foliation on a manifold $M$. Let
$\Diff(M,\FF)$ be the group of foliated diffeomorphisms
$(M,\FF)\to(M,\FF)$. The elements of $\Diff(M,\FF)$ that are
\cinf\ leafwisely homotopic to $\id_M$ form a normal subgroup
$\Diff_0(\FF)$, and let $\overline{\Diff}(M,\FF)$ denote the
corresponding quotient group. A {\em right transverse action\/} of
a group $G$ on $(M,\FF)$ is an anti-homomorphism
$\Phi:G\to\overline{\Diff}(M,\FF)$. A {\em local representation\/}
of $\Phi$ on some open subset $O\subset G$ is a map $\phi:M\times
O\to M$ such that $\phi_g=\phi(\cdot,g)$ is a foliated
diffeomorphism representing $\Phi_g$ for all $g\in G$. Then $\Phi$
is said to be of {\em class \cinf\/} if it has a \cinf\ local
representation on each small enough open subset of $G$.

Recall that the leafwise de Rham complex $(\Omega(\FF),d_\FF)$
consists of the differential forms on the leaves which are \cinf\
on $M$, endowed with the de~Rham derivative of the leaves. Its
cohomology $H(\FF)$ is called the leafwise cohomology. This
becomes a topological vector space with the topology induced by
the \cinf\ topology, and its maximal Hausdorff quotient  is the
reduced leafwise cohomology $\overline H(\FF)$.

Consider the canonical right action of  $\Diff(M,\FF)$ on
$\overline H(\FF)$ defined by pulling-back leafwise differential
forms. Since $\Diff_0(\FF)$ acts trivially, we get a canonical
right action of $\overline{\Diff}(M,\FF)$ on $\overline H(\FF)$.
Then any right transverse action $\Phi$ of a group $G$ on
$(M,\FF)$ induces a left action $\Phi^*$ of $G$ on $\overline
H(\FF)$.

Suppose from now on that $\FF$ is a Lie foliation and the manifold
$M$ is closed. It is shown that its transverse Lie structure can be
described as a right transverse action $\Phi$ of its structural Lie
group $G$ on $(M,\FF)$. Consider the induced left action $\Phi^*$ of
$G$ on $\overline H(\FF)$. For each $g\in G$, we would like to
define the supertrace $\Tr^s\Phi^*_g$, which could be called the
{\em leafwise Lefschetz number\/} $L(\Phi_g)$ of $\Phi_g$. This can
be achieved when $\overline H(\FF)$ is of finite dimension,
obtaining a \cinf\ function $L(\FF)$ on $G$ defined by
$L(\FF)(g)=L(\Phi_g)$; the value of $L(\FF)$ at the identity element
$e$ of $G$ is the Euler characteristic $\chi(\FF)$ of $\overline
H(\FF)$, which can be called the {\em leafwise Euler
characteristic\/} of $\FF$. But $\overline H(\FF)$ may be of
infinite dimension, even when the leaves are dense \cite{AlvHector},
and thus $L(\FF)$ is not defined in general.

The first goal of this paper is to show that, in general, the role
of the function $L(\FF)$ can be played by a distribution
$L_{\text{\rm dis}}(\FF)$ on $G$, called the {\em Lefschetz
distribution\/} of $\FF$, whose singularities are motivated by the
infinite dimension of $\overline H(\FF)$.

The first ingredient to define $L_{\text{\rm dis}}(\FF)$ is the
leafwise Hodge theory studied in \cite{AlvKordy:heat} for
Riemannian foliations; recall that Lie foliations form a specially
important class of Riemannian foliations \cite{Molino82}. Fix a
bundle-like metric on $M$ whose transverse part is induced by a
left invariant Riemannian metric on $G$. For the induced
Riemannian structure on the leaves, let $\D_\FF$ be the Laplacian
of the leaves operating in $\Omega(\FF)$. The kernel $\HH(\FF)$ of
$\D_F$ is the space of harmonic forms on the leaves that are
\cinf\ on $M$. The metric induces an $L^2$ inner product on
$\Omega(\FF)$, obtaining a Hilbert space $\bfOmega(\FF)$. Then
$\D_\FF$ is an essentially self-adjoint operator in
$\bfOmega(\FF)$ whose closure is denoted by $\bfDelta_\FF$. The
kernel of $\bfDelta_\FF$ is denoted by $\bfHH(\FF)$, and let
$\bfPi:\bfOmega(\FF)\to\bfHH(\FF)$ denote the orthogonal
projection. In \cite{AlvKordy:heat}, it is proved that $\bfPi$ has
a restriction $\Pi:\Omega(\FF)\to\HH(\FF)$ that induces an
isomorphism $\overline H(\FF)\cong\HH(\FF)$, which can be called
the {\em leafwise Hodge isomorphism\/}.

Let $\Lambda$ be the volume form of $G$, and let $\phi:M\times
O\to M$ be a \cinf\ local representation of $\Phi$. For each
$f\in\cinfc(O)$, consider the operator
$$
P_f=\int_G\phi_g^*\cdot f(g)\,\Lambda(g)\circ\Pi
$$
in $\Omega(\FF)$. Our first main result is the following.

\begin{prop}\label{p:trace class}
$P_f$ is of trace class, and the functional $f\mapsto\Tr^s P_f$
defines a distribution on $O$.
\end{prop}

It can be easily seen that $\Tr^s P_f$ is independent of the
choice of $\phi$, and thus the distributions given by
Proposition~\ref{p:trace class} can be combined to define a
distribution $L_{\text{\rm dis}}(\FF)$ on $G$; this is the {\em
Lefschetz distribution\/} of $\FF$.

Observe that $L_{\text{\rm dis}}(\FF)\equiv L(\FF)\cdot\Lambda$
when $\overline H(\FF)$ is of finite dimension. This justifies the
consideration of $L_{\text{\rm dis}}(\FF)$ as a generalization of
$L(\FF)$; in particular, the germ of $L_{\text{\rm dis}}(\FF)$ at
$e$ generalizes $\chi(\FF)$.

 If the operators $P_f$ are restricted to $\Omega^i(\FF)$ for each degree
$i$, its trace defines a distribution
$\Tr^i_{\text{\rm dis}}(\FF)$, called {\em distributional trace\/}, whose germ at $e$ generalizes the
{\em leafwise Betti number\/} $\beta^i(\FF)=\dim\overline
H^i(\FF)$.

The distributions $L_{\text{\rm dis}}(\FF)$ and $\Tr^i_{\text{\rm
dis}}(\FF)$ depend on $\Lambda$ and $\FF$, endowed with the
transverse Lie structure. If the leaves are dense, then the
transverse Lie structure is determined by the foliation, and thus
these distributions depend only on $\Lambda$ and the foliation. On
the other hand, the dependence on $\Lambda$ can be avoided by
using top dimensional currents instead of distributions, in the
obvious way.

Our second goal is to prove a distributional version of the
Gauss-Bonett theorem, which describes $L_{\text{\rm dis}}(\FF)$
around $e$. Let $R_\FF$ be the curvature of the leafwise metric.
Suppose for simplicity that $\FF$ is oriented. Then
$\Pf(R_\FF/2\pi)\in \Omega^p(\FF)$ ($p=\dim\FF$) can be called the
{\em leafwise Euler form\/}. This form can be paired with
$\Lambda$, considered as a transverse invariant measure, to give a
differential form $\omega_\Lambda\wedge\Pf(R_\FF/2\pi)$ of top
degree on $M$. In particular, if $\dim\FF=2$, then
$$
\omega_\Lambda\wedge\Pf(R_\FF/2\pi)=\frac{1}{2\pi}\,K_\FF\,\omega_M\;,
$$
where $K_\FF$ is the Gauss curvature of the leaves and $\omega_M$
is the volume form of $M$. Let $\delta_e$ denote the Dirac measure
at $e$.

\begin{thm}[Distributional Gauss-Bonett theorem]\label{t:Gauss-Bonett}
We have
$$
L_{\text{\rm
dis}}(\FF)=\int_M\omega_\Lambda\wedge\Pf(R_\FF/2\pi)\cdot\delta_e
$$
on some neighborhood of $e$.
\end{thm}

To prove Theorem~\ref{t:Gauss-Bonett}, we really prove that
\begin{equation}\label{e:L dis=chi Lambda}
L_{\text{\rm dis}}(\FF)=\chi_\Lambda(\FF)\cdot\delta_e
\end{equation}
around $e$, where $\Lambda$ is considered as a transverse
invariant measure of $\FF$, and $\chi_\Lambda(\FF)$ is the
$\Lambda$-Euler characteristic of $\FF$ introduced by Connes
\cite{Co79}. Then Theorem~\ref{t:Gauss-Bonett} follows from the
index theorem of \cite{Co79}.

The third goal is to prove a distributional version of the
Lefschetz trace formula, which describes $L_{\text{\rm dis}}(\FF)$
on any small enough open subset of $G$. For a \cinf\ local
representation $\phi:M\times O\to M$ of $\Phi$, let $\phi':M\times
O\to M\times O$ be the map defined by $\phi'(x,g)=(\phi_g(x),g)$.
The fixed point set of $\phi'$, $\Fix(\phi')$, consists of the
points $(x,g)$ such that $\phi_g(x)=x$. A point
$(x,g)\in\Fix(\phi')$ is said to be {\em leafwise simple\/} when
$\phi_{g*}-\id:T_x\FF\to T_x\FF$ is an isomorphism; in this case,
the sign of the determinant of this isomorphism is denoted by
$\epsilon(x,g)$. The set of leafwise simple fixed points of
$\phi'$ is denoted by $\Fix_0(\phi')$. Let $\pr_1:M\times O\to M$
and $\pr_2:M\times O\to O$ be the factor projections. It is proved
that $\Fix_0(\phi')$ is a \cinf\ manifold of dimension equal to
$\codim\FF$. Moreover the restriction $\pr_1:\Fix_0(\phi')\to M$
is a local embedding transverse to $\FF$. So $\Lambda$ defines a
measure $\Lambda'_{\Fix_0(\phi')}$ on $\Fix_0(\phi')$. Observe
that $\pr_2:\Fix(\phi')\to O$ is a proper map.

\begin{thm}[Distributional Lefschetz trace formula]\label{t:Lefschetz}
Suppose that every fixed point of $\phi'$ is leafwise simple. Then
$$
L_{\text{\rm
dis}}(\FF)=\pr_{2*}(\epsilon\cdot\Lambda'_{\Fix(\phi')})
$$
on $O$.
\end{thm}

To prove Theorem~\ref{t:Lefschetz}, we consider certain
submanifold $M'_1\subset M\times O$ endowed with a foliation
$\FF'_1$, whose leaves are of the form $L\times\{g\}$, where $L$
is a leaf of $\FF$ and $g\in G$. It is proved that $\pr_2(M'_1)$
is open in some orbit of the adjoint action of $G$ on itself,
$\pr_1:M'_1\to M$ is a local diffeomorphism, and
$\FF'_1=\pr_1^*\FF$. So $\Lambda$ lifts to a transverse invariant
measure $\Lambda'_1$ of $\FF'_1$. Moreover the restriction
$\phi'_1$ of $\phi'$ to $M'_1$ is defined and maps each leaf of
$\FF'_1$ to itself. For each $f\in\cinfc(O)$ supported in an
appropriate open subset $O_1\subset O$, the transverse invariant
measure $\Lambda'_{1,f}=\pr_2^*f\cdot \Lambda'_1$ is compactly
supported. Then the $\Lambda'_{1,f}$-Lefschetz number
$L_{\Lambda'_{1,f}}(\phi'_1)$ is defined according to
\cite{HeitschLazarov}. Without assuming any condition on the fixed
point set, we show that
\begin{equation}\label{e:L dis=L Lambda' 1,f}
\langle L_{\text{\rm
dis}}(\FF),f\rangle=L_{\Lambda'_{1,f}}(\phi'_1)\;.
\end{equation}
We have that $\Fix(\phi'_1)$ is a \cinf\ local transversal of
$\FF'_1$. Hence Theorem~\ref{t:Lefschetz} follows from~\eqref{e:L
dis=L Lambda' 1,f} and the foliation Lefschetz theorem of
\cite{HeitschLazarov,Rich}.

The numbers $\chi_\Lambda(\FF)$ and $L_{\Lambda'_{1,f}}(\phi'_1)$
are defined by using $L^2$ differential forms on the leaves,
whilst $L_{\text{\rm dis}}(\FF)$ is defined by using leafwise
differential forms that are \cinf\ on $M$. These are sharply
different conditions when the leaves are not compact.
So~\eqref{e:L dis=chi Lambda} and~\eqref{e:L dis=L Lambda' 1,f}
are surprising relations.

By~\eqref{e:L dis=L Lambda' 1,f}, $L_{\text{\rm dis}}(\FF)$ is
supported in the union of a discrete set of orbits of the adjoint
action. Therefore, when $\codim\FF>0$, $L_{\text{\rm dis}}(\FF)$
is \cinf\ just when it is trivial, obtaining the following.

\begin{cor}\label{c:L=0}
If $\overline H(\FF)$ is of finite dimension and $\codim\FF>0$,
then $L_{\text{\rm dis}}(\FF)\equiv L(\FF)=0$.
\end{cor}

By Corollary~\ref{c:L=0}, $\chi(\FF)$ is useless: it vanishes just
when it can be defined. Moreover $\chi_\Lambda(\FF)=0$ in this
case by~\eqref{e:L dis=chi Lambda}. So, when $\codim\FF>0$, the
condition $\chi_\Lambda(\FF)\neq0$ yields $\dim\overline
H(\FF)=\infty$. More precise results of this type would be
desirable.

Let $\dim\FF=p$. When the leaves are dense, $\beta^0(\FF)$ and
$\beta^p(\FF)$ are finite, and thus $\Tr_{\text{\rm dis}}^0(\FF)$
and $\Tr_{\text{\rm dis}}^p(\FF)$ are \cinf. On the other hand,
when the leaves are not compact, the $\Lambda$-Betti numbers of
\cite{Co79} satisfy $\beta_\Lambda^0(\FF)=\beta_\Lambda^p(\FF)=0$.
Then the following result follows from~\eqref{e:L dis=chi Lambda}
and Corollary~\ref{c:L=0}.

\begin{cor}\label{c:beta dis 1}
If $\codim\FF>0$, $\dim\FF=2$ and the leaves are dense, then
$\Tr_{\text{\rm dis}}^1(\FF)-\beta^1_\Lambda(\FF)\cdot\delta_e$ is
\cinf\ around $e$.
\end{cor}

In Corollary~\ref{c:beta dis 1}, we could say that
$\beta^1_\Lambda(\FF)\cdot\delta_e$ is the ``singular part'' of
$\Tr_{\text{\rm dis}}^1(\FF)$ around $e$.

\begin{cor}\label{c:dim overline H 1(FF)=infty}
Suppose that $\codim\FF>0$ and $\dim\FF=2$. If there is a
nontrivial harmonic $L^2$ differential form of degree one on some
leaf, then $\dim\overline H^1(\FF)=\infty$.
\end{cor}

It would be nice to generalize Corollary~\ref{c:dim overline H
1(FF)=infty} for arbitrary dimension. Thus we conjecture the
following.

\begin{conj}\label{q:beta dis i}
If $\codim\FF>0$ and the leaves are dense, then $\Tr_{\text{\rm
dis}}^i(\FF)-\beta^i_\Lambda(\FF)\cdot\delta_e$ is \cinf\ around $e$
for each degree $i$.
\end{conj}

The main results were proved in \cite{AlvKordy:Betti} for the case
of codimension one. Our results also overlap the corresponding
results of \cite{Mumken}.

We hope to prove elsewhere another version of
Theorem~\ref{t:Lefschetz} with a more general condition on the fixed
points, always satisfied by some local representation $\phi$ of
$\Phi$ defined around any point of $G$. By~\eqref{e:L dis=L Lambda'
1,f}, what is needed is another version of the Lefschetz theorem of
\cite{HeitschLazarov}, which holds for more general fixed point sets
when the transverse measure is \cinf.

The idea of using such type of trace class operators to define
distributional spectral invariants is due to Atiyah and Singer
\cite{Atiyah74,Singer:recent}. They consider transversally elliptic
operators with respect to compact Lie group actions. Further
generalizations to foliations and non-compact Lie group actions
were given in \cite{Nestke-Z,Co:nc,trans,noncom}. In our case,
$\D_\FF$ is not transversally elliptic with respect to any Lie
group action or any foliation, but it can be considered as being
``transversely elliptic'' with respect to the structural
transverse action; this simply means that it is elliptic along the
leaves of $\FF$.

\section{Transverse actions}\label{s:actions}

Recall that a foliation $\FF$ on a manifold $M$ can be described by
a {\em foliated cocycle\/}, which is a collection $\{U_i,f_i\}$,
where $\{U_i\}$ is an open cover of $X$ and each $f_i$ is a
topological submersion of $U_i$ onto some manifold $T_i$ whose
fibers are connected open subsets of $\R^n$, such that the following
{\em compatibility condition\/} is satisfied: for every $x\in
U_i\cap U_j$, there is an open neighborhood $U_{i,j}^x$ of $x$ in
$U_i\cap U_j$ and a homeomorphism $h_{i,j}^x:f_i(U_{i,j}^x)\to
f_j(U_{i,j}^x)$ such that $f_j=h_{i,j}^x\circ f_i$ on $U_{i,j}^x$.
Two foliated cocycles describe the same foliation $\FF$ when their
union is a foliated cocycle. The {\em leaf topology\/} on $M$ is the
topology with a base given by the open sets of the fibers of all the
submersions $f_i$. The {\em leaves\/} of $\FF$ are the connected
components of $M$ with the leaf topology. The leaf through each
point $x\in M$ is denoted by $L_x$. The pseudogroup on
$\bigsqcup_iT_i$ generated by the maps $h^x_{i,j}$, given by the
compatibility condition, is called (a representative of) the {\em
holonomy pseudogroup\/} of $\FF$, and describes the ``transverse
dynamics'' of $\FF$. Different foliated cocycles of $\FF$ induce
equivalent pseudogroups in the sense of \cite{Haefliger:80minimal,Haefliger:02compactly}.

Another representative of the
holonomy pseudogroup is defined on any transversal of $\FF$ that
meets every leaf. It is generated by ``sliding'' small open subsets
(local transversals) along the leaves; its precise
definition is given in \cite{Haefliger:80minimal}.

When $M$ is a \cinf\ manifold, it is said that $\FF$ is \cinf\ if
it is described by a foliated cocycle $\{U_i,f_i\}$ which is
\cinf\ in the sense that each $f_i$ is a \cinf\ submersion to some
\cinf\ manifold.

Let $\Gamma$ be a group of homeomorphisms of a manifold $T$. A
foliated cocycle $(U_i,f_i)$ of $\FF$, with $f_i:U_i\to T_i$, is
said to be {\em $(T,\Gamma)$-valued\/} when each $T_i$ is an open
subset of $T$, and the maps $h^x_{i,j}$, given by the
compatibility condition, are restrictions of maps in $\Gamma$. A
{\em transverse $(T,\Gamma)$-structure\/} of $\FF$ is given by a
$(T,\Gamma)$-valued foliated cocycle, and two $(T,\Gamma)$-valued
foliated cocycles define the same transverse
$(T,\Gamma)$-structure when their union is a $(T,\Gamma)$-valued
foliated cocycle. When $\FF$ is endowed with a transverse
$(T,\Gamma)$-structure, it is called a {\em
$(T,\Gamma)$-foliation\/}.

Let $\FF$ and $\GG$ be foliations on manifolds $M$ and $N$,
respectively. Recall the following concepts. A {\em foliated
map\/} $f:(M,\FF)\to(N,\GG)$ is a map $f:M\to N$ that maps each
leaf of $\FF$ to a leaf of $\GG$; the simpler notation $f:\FF\to
\GG$ will be also used. A {\em leafwise homotopy\/} (or {\em integrable homotopy\/}) between two
continuous foliated maps $f,f':(M,\FF)\to(N,\GG)$ is a continuous
map $H:M\times I\to N$ ($I=[0,1]$) such that the path
$H(x,\cdot):I\to N$ lies in a leaf of $\GG$ for each $x\in M$; in
this case, it is said that $f$ and $f'$ are {\em leafwisely
homotopic\/} (or {\em integrably homotopic\/}).

Suppose from now on that $\FF$ and $\GG$ are \cinf. Two \cinf\
foliated maps are said to be {\em \cinf\ leafwisely homotopic\/}
when there is a \cinf\ leafwise homotopy between them. As usual,
$T\FF\subset TM$ denotes the subbundle of vectors tangent to the
leaves of $\FF$, $\mathfrak{X}(M,\FF)$ denotes the Lie algebra of
infinitesimal transformations of $(M,\FF)$, and
$\mathfrak{X}(\FF)\subset\mathfrak{X}(M,\FF)$ is the normal Lie
subalgebra of vector fields tangent to the leaves of $\FF$ (\cinf\
sections of $T\FF\to M$). Then we can consider the quotient Lie
algebra
$\overline{\mathfrak{X}}(M,\FF)=\mathfrak{X}(M,\FF)/\mathfrak{X}(\FF)$,
whose elements are called {\em transverse vector fields}. Observe
that, for each $x\in M$, the evaluation map
$\ev_x:\mathfrak{X}(M,\FF)\to T_xM$ induces a map
$\overline{\ev}_x:\overline{\mathfrak{X}}(M,\FF)\to T_xM/T_x\FF$,
which can be also called {\em evaluation map\/}. For any Lie algebra
$\mathfrak{g}$, a homomorphism
$\mathfrak{g}\to\overline{\mathfrak{X}}(M,\FF)$ is called an {\em
infinitesimal transverse action\/} of $\mathfrak{g}$ on $(M,\FF)$.
In particular, we have a canonical infinitesimal transverse action
of $\overline{\mathfrak{X}}(M,\FF)$ on $(M,\FF)$.

Let $\Diff(M,\FF)$ be the group of \cinf\ foliated diffeomorphisms
$(M,\FF)\to(M,\FF)$ with the operation of composition, let
$\Diff(\FF)\subset\Diff(M,\FF)$ be the normal subgroup \cinf\
foliated diffeomorphisms that preserve each leaf of $\FF$, and let
$\Diff_0(\FF)\subset\Diff(\FF)$ be the normal subgroup of \cinf\
foliated diffeomorphisms that are \cinf\ leafwisely homotopic to the
identity map. Then we can consider the quotient group
$\overline{\Diff}(M,\FF)=\Diff(M,\FF)/\Diff_0(\FF)$, whose operation
is also denoted by ``$\circ$''. The elements of
$\overline{\Diff}(M,\FF)$ can be called {\em transverse
transformations\/} of $(M,\FF)$. For any group $G$, an
anti-homomorphism $\Phi:G\to\overline{\Diff}(M,\FF)$,
$g\mapsto\Phi_g$, is called a {\em right transverse action\/} of $G$
on $(M,\FF)$. For an open subset $O\subset G$, a map $\phi:M\times
O\to M$ is called a {\em local representation\/} of $\Phi$ on $O$ if
$\phi_g=\phi(\cdot,g)\in\Phi_g$ for all $g\in O$. For any leaf $L$
of $\FF$ and any $g\in O$, the leaf $\phi_g(L)$ is independent of
the local representative $\phi$, and thus it will be denoted by
$\Phi_g(L)$. When $G$ is a Lie group, $\Phi$ is said to be of {\em
class\/} \cinf\ if it has a \cinf\ local representation around each
element of $G$.

Somehow, we can think of $\overline{\Diff}(M,\FF)$ as a Lie group
whose Lie algebra is $\overline{\mathfrak{X}}(M,\FF)$; indeed, it
will be proved elsewhere that, if $G$ is a simply connected Lie
group and $\mathfrak{g}$ is its Lie algebra of left invariant
vector fields, then there is a canonical bijection between
infinitesimal transverse actions of $\mathfrak{g}$ on $(M,\FF)$
and \cinf\ right transverse actions of $G$ on $(M,\FF)$.

The {\em leafwise de~Rham complex\/} $(\Omega(\FF),d_\FF)$ is the
space of differential forms on the leaves smooth on $M$ (\cinf\
sections of $\bigwedge T\FF^*\to M$) endowed with  the leafwise de
Rham differential. It is also a topological vector space with the
\cinf\ topology, and $d_\FF$ is continuous. The cohomology
$H(\FF)$ of $(\Omega(\FF),d_\FF)$ is called the {\em leafwise
cohomology\/} of $\FF$, which is a topological vector space with
the induced topology. Its maximal Hausdorff quotient
$\overline{H}(\FF)=H(\FF)/\overline{0}$ is called the {\em reduced
leafwise cohomology\/}.

By pulling back leafwise differential forms, any $\cinf$ foliated
map $f:(M,\FF)\to(N,\GG)$ induces a continuous homomorphism of
complexes, $f^*:\Omega(\GG)\to\Omega(\FF)$, obtaining a continuous
homomorphism $f^*:\overline{H}(\GG)\to\overline{H}(\FF)$.
Moreover, if $f$ is $\cinf$ leafwisely homotopic to another
$\cinf$ foliated map $f':(M,\FF)\to(M,\FF)$, then
$f^*={f'}^*:\overline H(\GG)\to\overline H(\FF)$ by standard
arguments \cite{BottTu}. Therefore, for any
$F\in\overline{\Diff}(M,\FF)$ and any $f\in F$, the endomorphism
$f^*$ of $\overline{H}(\FF)$ can be denoted by $F^*$. So any right
transverse action $\Phi$ of a group $G$ on $(M,\FF)$ induces a
left action $\Phi^*$ of $G$ on $\overline{H}(\FF)$ given by
$(g,\xi)\mapsto \Phi_g^*\xi$.

\section{Lie foliations}\label{s:Lie}
Let $\FF$ be a \cinf\ foliation of codimension $q$ on a \cinf\
closed manifold $M$. Let $G$ be a simply connected Lie group of
dimension $q$, and $\mathfrak{g}$ its Lie algebra of left
invariant vector fields. A {\em transverse Lie structure\/} of
$\FF$, with {\em structural Lie group\/} $G$ and {\em structural
Lie algebra\/} $\mathfrak{g}$, can be described with any of the
following objects that determine each other
\cite{Fedida,Molino82}:
\begin{itemize}

\item[(L.1)] A transverse $(G,G)$-structure of $\FF$, where
$G$ is identified with the group of its left translations.

\item[(L.2)] A $\mathfrak g$-valued $1$-form $\omega$ on $M$ such that
$\omega_x:T_xM\rightarrow {\mathfrak g}$ is surjective with kernel
$T_x\FF$ for every $x\in M$, and
$$
d\omega+\frac{1}{2}\,[\omega,\omega]=0\;.
$$

\item[(L.3)]
A homomorphism
$\theta:\mathfrak{g}\to\overline{\mathfrak{X}}(M,\FF)$ such that
the composite
$$
\begin{CD}
\mathfrak{g} @>{\theta}>> \overline{\mathfrak{X}}(M,\FF)
@>{\overline{\ev}_x}>> T_xM/T_x\FF
\end{CD}
$$
is an isomorphism for every $x\in M$.

\end{itemize}
In~(L.1), the elements of $G$ whose corresponding left translations
are involved in the definition of the transverse $(G,G)$-structure
form a subgroup $\Gamma$, which is called the {\em holonomy group\/}
of $\FF$. So the transverse $(G,G)$-structure is a transverse
$(G,\Gamma)$-structure. In~(L.2) and~(L.3), $\omega$ and $\theta$
can be respectively called the {\em structural form\/} and the {\em
structural infinitesimal transverse action\/}.

A \cinf\ foliation endowed with a transverse Lie structure is
called a {\em Lie foliation\/}; the terms {\em Lie
$G$-foliation\/} or {\em Lie $\mathfrak{g}$-foliation\/} are used
too. If the leaves are dense, then the transverse Lie structure is
unique, and thus it is determined by the foliation.

A Lie $G$-foliation $\FF$ on a \cinf\ closed manifold $M$ has the
following description due to Fedida \cite{Fedida,Molino82}. There
exists a regular covering $\pi :\widetilde{M} \rightarrow M$, a
fibre bundle $D:\widetilde{M}\rightarrow G$ and an injective
homomorphism $h:\Aut(\pi)\to G$ such that the leaves of
$\widetilde{\FF}=\pi^*\FF$ are the fibres of $D$, and $D$ is
$h$-equivariant; {\it i.e.\/},
$$
D\circ\sigma(\tilde{x}) =h(\sigma)\cdot D(\tilde{x})
$$
for all $\tilde{x}\in\widetilde{M}$ and
$\sigma\in\operatorname{Aut}(\pi)$. This $h$ is called the {\em
holonomy homomorphism\/}. By using the covering space
$\ker(h)\backslash\widetilde{M}$ of $M$ if necessary, we can
assume that $h$ is injective, and thus $\pi$ restricts to
diffeomorphisms of the leaves of $\widetilde{\FF}$ to the leaves
of $\FF$. The leaf of $\widetilde\FF$ through each point $\tilde
x\in\widetilde M$ will be denoted by $\widetilde L_{\tilde x}$.

Given a $(G,G)$-valued foliated cocycle $\{U_i,f_i\}$ defining the
transverse Lie structure according to~(L.1), the
$\mathfrak{g}$-valued $1$-form $\omega$ of~(L.2) and the
infinitesimal transverse action $\theta$ of~(L.3) can be defined as
follows. For $x\in U_i$ and $v\in T_xM$, $\omega_x(v)$ is the left
invariant vector field on $G$ whose value at $f_i(x)$ is
$f_{i*}(v)$. To define $\theta$, fix an auxiliary vector subbundle
$\nu\subset TM$ complementary of $T\FF$ ($TM=\nu\oplus T\FF$). Each
$X\in\mathfrak{g}$ defines a \cinf\ vector field
$X^\nu\in\mathfrak{X}(M,\FF)$ by the conditions $X^\nu(x)\in\nu_x$
and $f_{i*}(X^\nu(x))=X(f_i(x))$ if $x\in U_i$. Then $\theta(X)$ is
the class of $X^\nu$ in $\overline{\mathfrak{X}}(M,\FF)$, which is
independent of the choice of $\nu$.

By using Fedida's geometric description of $\FF$, the definitions
of $\omega$ and $X^\nu$ can be better understood:
\begin{itemize}

\item Let $\omega_G$ be the canonical $\mathfrak{g}$-valued $1$-form
on $G$ defined by $\omega_G(X(g))=X$ for any $X\in\mathfrak{g}$ and
any $g\in G$. Then $\omega$ is determined by the condition $\pi^*\omega=D^*\omega_G$.

\item Let $\tilde\nu=\pi_*^{-1}(\nu)\subset T\widetilde{M}$, which is
a vector subbundle complementary of $T\widetilde{\FF}$. Then, for any
$X\in\mathfrak{g}$, there is a unique $\widetilde{X}^\nu\in
\mathfrak{X}(\widetilde{M},\widetilde{\FF})$ which is a section of
$\tilde\nu$ and satisfies $D_*\circ\widetilde{X}^\nu=X\circ D$.
Since $D$ is $h$-equivariant, $\widetilde{X}^\nu$ is
$\Aut(\pi)$-invariant. Then $X^\nu$ is the projection of
$\widetilde{X}^\nu$ to $M$.

\end{itemize}

\section{Structural transverse action}\label{s:structural action}

Let $G$ be a simply connected Lie group, and let $\FF$ be a Lie
$G$-foliation on a closed manifold $M$. According to
Section~\ref{s:actions}, the structural infinitesimal transverse
action corresponds to a unique right transverse action of $G$ on
$(M,\FF)$, obtaining another description of the transverse Lie
structure:
\begin{itemize}

\item[(L.4)] A \cinf\ right transverse action $\Phi$ of $G$ on $(M,\FF)$
which has a  \cinf\ local representation $\phi$ around the identity element
$e$ of $G$ such that the composite
$$
\begin{CD}
T_eG @>{\phi^x_*}>> T_xM @>>> T_xM/T_x\FF
\end{CD}
$$
is an isomorphism for all $x\in M$, where $\phi^x=\phi(x,\cdot)$
and the second map is the canonical projection. This condition is
independent of the choice of $\phi$. This $\Phi$ is called the
{\em structural transverse action\/}.

\end{itemize}

To describe $\Phi$, consider Fedida's geometric description of
$\FF$ (Section~\ref{s:Lie}). For any $g\in G$, take a continuous,
piecewise \cinf\ path $c:I\rightarrow G$ with $c(0)=e$ and
$c(1)=g$. For any $\tilde{x}\in\widetilde{M}$, there exists a
unique continuous piecewise \cinf\ path
$\tilde{c}_{\tilde{x}}^{\nu}:I\rightarrow \widetilde{M}$ such that
\begin{itemize}

\item $\tilde{c}_{\tilde{x}}^{\nu}(0)=\tilde{x}$,

\item $\tilde{c}_{\tilde{x}}^{\nu}$ is tangent to $\tilde{\nu}$ at every
$t\in I$ where it is \cinf, and

\item $D\circ \tilde{c}_{\tilde{x}}^{\nu}(t)=D(\tilde{x})\cdot c(t)$ for any $t\in I$.

\end{itemize}
It is easy to see that such a $\tilde{c}_{\tilde{x}}^{\nu}$
depends smoothly on $\tilde{x}$.

\begin{lem}\label{l:tilde c}
We have $\sigma\circ
\tilde{c}_{\tilde{x}}^{\nu}=\tilde{c}_{\sigma(\tilde{x})}^{\nu}$
for $\tilde{x}\in\widetilde{M}$ and $\sigma\in
\operatorname{Aut}(\pi)$.
\end{lem}

\begin{proof}
This is a direct consequence of the $h$-equivariance of $D$ and
the unicity of the paths $\tilde{c}_{\tilde{x}}^{\nu}$.
\end{proof}

For each $g\in G$, let
$\tilde{\phi}_g:(\widetilde{M},\widetilde{\FF})\rightarrow(\widetilde{M},\widetilde{\FF})$
be the \cinf\ foliated diffeomorphism given by
$\tilde{\phi}_g(\tilde{x})=\tilde{c}_{\tilde{x}}^{\nu}(1)$. For
any $\tilde{x}\in\widetilde{M}$ and $\sigma\in
\operatorname{Aut}(\pi)$, we have
$$
\sigma\circ\tilde{\phi}_g(\tilde{x})=\sigma\circ\tilde{c}_{\tilde{x}}^\nu(1)=
\tilde{c}_{\sigma(\tilde{x})}^\nu(1)=\tilde{\phi}_g\circ
\sigma(\tilde{x})
$$
by Lemma~\ref{l:tilde c}, yielding $\sigma\circ
\tilde{\phi}_g=\tilde{\phi}_g\circ \sigma$. Therefore, there
exists a unique \cinf\ foliated diffeomorphism
$\phi_g:(M,\FF)\rightarrow (M,\FF)$ such that
$\pi\circ\tilde{\phi}_g=\phi_g\circ\pi$.

\begin{lem}\label{l:indepent of c}
The \cinf\ leafwise homotopy class of $\phi_g$ is independent of
the choice of $c$.
\end{lem}

\begin{proof}
Let $d:I\to G$ be another continuous and piecewise smooth path
with $d(0)=e$ and $d(1)=g$, which defines a \cinf\ foliated map
$\varphi_g:(M,\FF)\to(M,\FF)$ as above. Since $G$ is simply
connected, there exists a family of continuous and piecewise
smooth paths $c_s:I\to G$, depending smoothly on $s\in I$, with
$c_s(0)=e$, $c_s(1)=g$, $c_0=c$ and $c_1=d$. The paths $c_s$
induce a family of \cinf\ foliated maps
$\phi_{g,s}:(M,\FF)\to(M,\FF)$ as above, defining a \cinf\
leafwise homotopy between $\phi_g$ and $\varphi_g$.
\end{proof}

\begin{lem}\label{l:indepent of nu}
The  $\cinf$ leafwise homotopy class of $\phi_g$ is independent of
the choice of $\nu$.
\end{lem}

\begin{proof}
Let $\nu'\subset TM$ be another vector subbundle complementary of
$T\FF$, which can be used to define a \cinf\ foliated map
$\phi'_g$ as above. It is easy to find a $\cinf$ deformation of
vector subbundles of $\nu_s\subset TM$ complementary of $T\FF$,
$s\in I$, with $\nu_0=\nu$ and $\nu_1=\nu'$. Then the foliated
maps $\phi_{g,s}$, induced by the vector bundles $\nu_s$ as above,
define a \cinf\ leafwise homotopy between $\phi_g$ and
$\phi^\prime_g$.
\end{proof}

Therefore, for each $g$, the \cinf\ leafwise homotopy class
$\Phi_g$ of $\phi_g$ depends only on $g$, $\FF$ and its transverse
Lie structure. So a map $\Phi:G\to\overline{\Diff}(M,\FF)$ is
given by $g\mapsto\Phi_g$.

\begin{lem}\label{l:right transv action}
$\Phi$ is a right transverse action of $G$ in $(M,\FF)$.
\end{lem}

\begin{proof}
Given $g_1, g_2\in G$, let $c_1,c_2:I\to G$ be continuous,
piecewise smooth paths such that $c_1(0)=c_2(0)=e$, $c_1(1)=g_1$
and $c_2(1)=g_2$, which are used to define $\phi_{g_1}$ and
$\phi_{g_2}$ as above. Let $c:I\to G$ be the path product of $c_1$
and $L_{g_1}\circ c_2$, where $L_{g_1}$ denotes the left
translation by $g_1$. We have $c(0)=e$ and $c(1)=g_1g_2$. We can
use this $c$ to define $\phi_{g_1g_2}$, obtaining
$\phi_{g_1g_2}=\phi_{g_2}\circ \phi_{g_1}$, and thus
$\Phi_{g_1g_2}=\Phi_{g_2}\circ \Phi_{g_1}$.
\end{proof}

\begin{lem}\label{l:cinf}
$\Phi$ is $\cinf$.
\end{lem}

\begin{proof}
It is easy to prove that each element of $G$ has a neighbourhood $O$
such that there is a \cinf\ map $c:I\times O\to G$ so that each
$c_g=c(\cdot,g)$ is a path from $e$ to $g$. The corresponding
foliated diffeomorphisms $\phi_g$ form a \cinf\ representation of
$\Phi$ on $O$.
\end{proof}

This construction defines the structural transverse action $\Phi$.
According to Section~\ref{s:actions}, $\Phi$ induces a left action
$\Phi^*$ of $G$ on $\overline{H}(\FF)$.

\begin{lem}\label{l:phi e=id M}
There is a local representation $\varphi:M\times O\to M$ of $\Phi$
around the identity element $e$ such that $\varphi_e=\id_M$.
\end{lem}

\begin{proof}
Construct $\phi$ like in the proof of Lemma~\ref{l:cinf} such that
$e\in O$ and $c_e$ is the constant path at $e$.
\end{proof}

Let $\varphi:M\times O\to M$ be a local representation of $\Phi$.
A map $\tilde\varphi:\widetilde M\times O\to\widetilde M$ is
called a {\em lift\/} of $\varphi$ if
$\pi\circ\tilde\varphi_g=\varphi_g\circ\pi$ for all $g\in O$,
where $\tilde\varphi_g=\tilde\varphi(\cdot,g)$. In particular, the
above construction of $\phi$ also gives a lift $\tilde\phi$. Let
$R_g:G\to G$ denote the right translation by any $g\in G$.

\begin{lem}\label{l:tilde varphi}
Any \cinf\ lift $\tilde\varphi:\widetilde M\times O\to\widetilde
M$ of each \cinf\ local representation $\varphi:M\times O\to M$ of
$\Phi$, such that $O$ is connected, satisfies
$D\circ\tilde\varphi_g=R_g\circ D$ for all $g\in O$.
\end{lem}

\begin{proof}
It is enough to prove the result when $O$ is as small as desired.
It is clear that the property of the statement is satisfied by the
maps $\tilde\phi$ constructed above for connected $O$.

For an arbitrary $\varphi$, if $O$ is small enough and connected,
there is some $\phi:M\times O\to M$ defined by the above
construction and some homotopy $H:M\times O\times I\to M$ between
$\varphi$ and $\phi$ such that each path $t\mapsto H(x,g,t)$ is
contained in a leaf of $\FF$. This $H$ lifts to a homotopy
$\widetilde H:\widetilde M\times O\times I\to\widetilde M$ between
$\tilde\varphi$ and $\tilde\phi$ so that each path
$t\mapsto\widetilde H(\tilde x,g,t)$ is contained in a leaf of
$\widetilde\FF$. Then $D\circ\tilde\varphi=D\circ\tilde\phi$,
completing the proof.
\end{proof}

\begin{cor}\label{c:tilde varphi}
$\tilde\varphi:\widetilde L\times O\to\widetilde M$ is a \cinf\
embedding for each leaf $\widetilde L$ of $\widetilde\FF$.
\end{cor}

The transverse Lie structure of $\FF$ lifts to a transverse Lie
structure of $\widetilde\FF$, whose structural right transverse
action is locally represented by the \cinf\ lifts of \cinf\ local
representations of $\Phi$.

\section{The Hodge isomorphism}\label{s:Hodge}

Recall that any Lie foliation is Riemannian \cite{Reinhart59}. Then
fix a bundle-like metric on $M$ \cite{Reinhart59}, and equip the
leaves of $\FF$ with the induced Riemannian metric. Let $\delta_\FF$
denote the leafwise coderivative on the leaves operating in
$\Omega(\FF)$, and set $D_\FF=d_\FF+\delta_\FF$. Then
$\Delta_\FF=D_\FF^2=d_\FF\circ d_\FF+d_\FF\circ\delta_\FF$ is the
leafwise Laplacian operating in $\Omega(\FF)$. Let
$\HH(\FF)=\ker\Delta_\FF$ (the space of leafwise harmonic forms
which are smooth on $M$). Since the metric is bundle-like, the
transverse volume element is holonomy invariant, which implies that
$D_\FF$ and $\Delta_\FF$ are symmetric, and thus they have the same
kernel.

Let $\bfOmega(\FF)$ be the Hilbert space of square integrable
leafwise differential forms on $M$. The metric of $M$ induces a
Hilbert structure in $\bfOmega(\FF)$. For any $\cinf$ foliated map
$f:(M,\FF)\to(M,\FF)$, the endomorphism $f^*$ of $\Omega(\FF)$ is
obviously $L^2$-bounded, and thus extends to a bounded operator
$f^*$ in $\bfOmega(\FF)$. Consider $D_\FF$ and $\Delta_\FF$ as
unbounded operators in $\bfOmega(\FF)$, which are essentially
self-adjoint \cite{Chernoff}, and whose closures are denoted by
$\bfD_\FF$ and $\bfDelta_\FF$ (see {\it e.g.\/} \cite{Alv-T,tang}). By
\cite{AlvKordy:heat}, $\bfHH(\FF)=\ker \bfDelta_\FF$ is the closure
of $\HH(\FF)$ in $\bfOmega(\FF)$, and the orthogonal projection
$\bfPi:\bfOmega(\FF) \to\bfHH(\FF)$ has a restriction $\Pi:
\Omega(\FF)\to\HH(\FF)$, which induces a leafwise Hodge isomorphism
$$
\overline{H}(\FF)\cong \HH(\FF)\;.
$$
For any \cinf\ foliated map $f:(M,\FF)\to(M,\FF)$, the
homomorphism $f^*:\overline{H}(\FF)\to\overline{H}(\FF)$
corresponds to the operator $\Pi\circ f^*$ in $\HH(\FF)$ via the
Hodge isomorphism. So the left $G$-action on $\overline{H}(\FF)$,
defined in Section~\ref{s:structural action}, corresponds to the
left $G$-action on $\HH(\FF)$ given by $(g,\alpha)\mapsto\Pi\circ
\phi^*_g\alpha$ for any $\phi_g\in\Phi_g$.

Since the left action of $G$ on $\HH(\FF)$ is $L^2$-continuous, we
get an extended left action of $G$ on $\bfHH(\FF)$ given by
$(g,\alpha)\mapsto\bfPi\circ \phi^*_g\alpha$ for any
$\phi_g\in\Phi_g$.

These actions on $\HH(\FF)$ and $\bfHH(\FF)$ are continuous on $G$
since $\Phi$ is \cinf.

\section{A class of smoothing operators}

\subsection{Preliminaries on smoothing and trace class operators}

Let $\omega_M$ denote the volume forms of $M$. A {\em smoothing
operator\/} in $\Omega(\FF)$ is a linear map
$P:\Omega(\FF)\to\Omega(\FF)$, continuous with respect to the
$\cinf$ topology, given by
$$
(P\alpha)(x)=\int_Mk(x,y)\,\alpha(y)\,\omega_M(y)
$$
for some \cinf\ section $k$ of $\bigwedge T\FF^*\boxtimes
\bigwedge T\FF$ over $M\times M$; thus
$$
k(x,y)\in \bigwedge T\FF^*_x\otimes \bigwedge T\FF_y \equiv
\operatorname{Hom}(\bigwedge T\FF^*_y,\bigwedge T\FF^*_x)
$$
 for any $x,y\in M$. This $k$ is called the {\em smoothing kernel\/}
 or {\em Schwartz kernel\/} of $P$. Such a $P$ defines a trace class
 operator in $\bfOmega(\FF)$, and we have
$$
\Tr P=\int_M \Tr k(x,x)\,\omega_M(x)\;.
$$
The supertrace formalism will be also used. For any homogeneous
operator $T$ in $\bfOmega(\FF)$ or in $\bigwedge T_x\FF^*$, let
$T^\pm$ denote its restriction to the even and odd degree part, and
let $T^{(i)}$ denote its restriction to the part of degree $i$. If
$T$ is of trace class, then its supertrace is
$$
\Tr^sT=\Tr
T^+-\Tr T^-=\sum_i(-1)^i\,\Tr T^{(i)}\;.
$$
Thus
$$
\Tr^sP=\int_M \Tr^s k(x,x)\,\omega_M(x)\;.
$$

Let $W^k\Omega(\FF)$ denote the Sobolev space of order $k$ of
leafwise differential forms on $M$, and let $\|\cdot\|_k$ denote a
norm of $W^k\Omega(\FF)$. A continuous operator $P$ in
$\Omega(\FF)$ is smoothing if and only if $P$ extends to a bounded
operator $P:W^k\Omega(\FF)\to W^l\Omega(\FF)$ for any $k$ and $l$.

If an operator $P$ in $\Omega(\FF)$ has an extension
$P:W^k\Omega(\FF)\to W^\ell\Omega(\FF)$, then $\|P\|_{k,\ell}$
denotes the norm of this extension; the notation $\|P\|_k$ is used
when $k=\ell$. By the Sobolev embedding theorem, the trace of a
smoothing operator $P$ in $\Omega(\FF)$ can be estimated in the
following way: for any $k>\dim M$, there is some $C>0$ independent
of $P$ such that
\begin{equation}\label{e:Tr}
|\Tr P|\leq C\,\|P\|_{0,k}\;.
\end{equation}

\subsection{The class $\DD$}\label{ss:DD}

Let $\AAA$ be the set of all functions $\psi:\R\to\C$, extending
to an entire function $\psi$ on $\C$ such that, for each compact
set $K\subset\R$, the set of functions $\{(x\mapsto\psi(x+iy)) \
|\  y\in K\}$ is bounded in the Schwartz space $\SSS(\R)$. This
$\AAA$ has a structure of Frechet algebra, and, in fact, it is a
module over $\C[z]$. This algebra contains all functions with
compactly supported Fourier transform, and the functions $x\mapsto
e^{-tx^2}$ with $t>0$.

By \cite[Proposition~4.1]{Roe87}, there exists a ``functional
calculus map'' $\AAA\to \End(\Omega(\FF))$,
$\psi\mapsto\psi(D_{\FF})$, which is a continuous homomorphism of
$\C[z]$-modules and of algebras. Any operator $\psi(D_{\FF})$,
$\psi\in\AAA$, extends to a bounded operator in $W^k\Omega(\FF)$
for any $k$ with the following estimate for its norm: there is
some $C>0$, independent of $\psi$, such that
\begin{equation}\label{e:fdnorm}
\|\psi(D_{\FF})\|_k\leq
\int|\hat{\psi}(\xi)|\,e^{C\,|\xi|}\,d\xi\;,
\end{equation}
where $\hat{\psi}$ denotes the Fourier transform of $\psi$.
Therefore, for any natural $N$, the operator
$(\id+\Delta_{\FF})^N\psi(D_{\FF})$ extends to a bounded operator
in $W^k\Omega(\FF)$ for any $k$ whose norm can be estimated as
follows: there is some $C>0$, independent of $\psi$, such that
\begin{equation}\label{e:fdnorm1}
\|(\id+\Delta_{\FF})^N\psi(D_{\FF})\|_k \leq\int
|(\id-\partial^2_{\xi})^N\hat{\psi}(\xi)|\,e^{C\,|\xi|}\,d\xi\;.
\end{equation}

Fix a left-invariant Riemannian metric on $G$, and let $\Lambda$
denote its volume form. We can assume that the metrics on $M$ and
$G$ agree in the sense that the maps $f_i$ of (L.1) are Riemannian
submersions (Section~\ref{s:Lie}). Thus $D:\widetilde M\to G$ is a
Riemannian submersion with respect to the lift of the bundle-like
metric to $\widetilde M$.

A {\em leafwise differential\/} operator in $\Omega(\FF)$ is a
differential operator which involves only leafwise derivatives; for
instance, $d_\FF$, $\delta_\FF$, $D_\FF$ and $\D_\FF$ are leafwise
differential operators. A family of leafwise differential operators
in $\Omega(\FF)$, $A=\{A_v\ |\ v\in V\}$, is said to be {\em
smooth\/} when $V$ is a \cinf\ manifold and, with respect to \cinf\
local coordinates, the local coefficients of each $A_v$ depend
smoothly on $v$ in the $C^\infty$-topology.  We also say that $A$ is
{\em compactly supported\/} when there is some compact
subset $K\subset V$ such that $A_v=0$ if $v\notin K$.
Given another smooth family of leafwise differential operators in
$\Omega(\FF)$ with the same parameter manifold, $B=\{B_v\ |\ v\in
V\}$, the {\em composite\/} $A\circ B$ is the family defined by
$(A\circ B)_v=A_v\circ B_v$. Similarly, we can define the {\em
sum\/} $A+B$ and the {\em product\/} $\lambda\cdot A$ for some
$\lambda\in\R$.

We introduce the class $\DD$ of operators
$P:\Omega(\FF)\to\Omega(\FF)$ of the form
$$
P=\int_O\phi^*_g\circ A_g\,\Lambda(g)\circ\psi(D_{\FF})\;,
$$
where $O$ is some open subset of $G$, $\phi:M\times O\to M$ is a
\cinf\ local representation of $\Phi$, $A=\{A_g \ |\  g\in O\}$ is
a smooth compactly supported family of leafwise differential
operators in $\Omega(\FF)$, and $\psi\in\AAA$.

\begin{prop}\label{p:smoothing}
Any operator $P\in \DD$ is a smoothing operator in $\Omega(\FF)$.
\end{prop}

\begin{proof}
Let $P\in \DD$ as above. By~Ê\eqref{e:fdnorm1} and since the
operator $\phi_g^*$ preserves any Sobolev space, $P$ defines a
bounded operator in $W^k\Omega(\FF)$ for any $k$.

Let $\varphi:M\times O_0\to M$ be a \cinf\ local representation of
$\Phi$ on some open neighborhood $O_0$ of the identity element
$e$; we can assume that $\varphi_e=\id_M$ by
Corollary~\ref{c:tilde varphi}. For any $Y\in \mathfrak g$, let
$\widehat{Y}$ be the first order differential operator in
$\Omega(\FF)$ defined by
$$
\widehat{Y}u=\left.\frac{d}{dt}\,\varphi^*_{\exp
tY}u\right|_{t=0}\;,
$$
which makes sense because $\exp tY\in O_0$ for any $t>0$ small
enough.

Fix a base $Y_1,\ldots,Y_q$ of $\mathfrak{g}$. Then the second
order differential operator $L=-\sum_{j=1}^q\widehat{Y}_j^2$ in
$\Omega(\FF)$ is transversely elliptic. Moreover $\D_\FF$ is
leafwise elliptic. By the elliptic regularity theorem, it suffices
to prove that $L^N\circ P$ and $\D_\FF^N\circ P$ belong to $\DD$
for any natural $N$. In turn, this follows by showing that $Q\circ
P$ and $\widehat{Y}\circ P$ are in $\DD$ for any leafwise
differential operator $Q$ and any $Y\in\mathfrak{g}$.

We have
$$
Q\circ P=\int_O\phi^*_g\circ  B_g\,\Lambda(g)\circ\psi(D_{\FF})\;,
$$
where $B_g=(\phi^*_g)^{-1}\circ Q\circ \phi^*_g\circ A_g$. Since
$\phi_g$ is a foliated map, it follows that $\{B_g \ |\  g\in O\}$
is a smooth family of leafwise differential operators, yielding
$Q\circ P\in\DD$.

For $g\in O$ and $a\in O_0$ close enough to $e$, let
$$
F_{a,g}=\phi_{ag}\circ\varphi_a\circ\phi_g^{-1}\;.
$$
Observe that $F_{e,g}=\id_M$ because $\varphi_e=\id_M$. For each
$Y\in\mathfrak{g}$, we get a smooth family $V_Y=\{V_{Y,g}\ |\ g\in
O\}$ of first order leafwise differential operators in
$\Omega(\FF)$ given by
$$
V_{Y,g}u=\left.\frac{d}{dt}\,F^*_{\exp tY,g}u\right|_{t=0}\;.
$$
Let also $L_YA=\{(L_YA)_g\ |\ g\in O\}$ be the smooth family of
leafwise differential operators given by
$$
(L_YA)_gu=\left.\frac{d}{dt}\,A_{\exp(-tY)\cdot g}u\right|_{t=0}.
$$
In particular, if $A_g$ is given by multiplication by $f(g)$ for
some $f\in\cinfc(G)$, then $(L_YA)_g$ is given by multiplication
by $(Yf)(g)$.

We proceed as follows:
\begin{align*}
\int_O\varphi^*_{\exp tY}\circ \phi^*_g\circ A_g\,\Lambda(g)
&=\int_O\phi^*_{\exp tY\cdot g}\circ F^*_{\exp tY,\exp(-tY)\cdot g}\circ A_g\,\Lambda(g)\\
&=\int_O\phi^*_{g}\circ F^*_{\exp tY,g}\circ A_{\exp tY\cdot
g}\,\Lambda(g)\;,
\end{align*}
yielding
\begin{align*}
\widehat{Y}\circ P
&=\lim_{t\to0}\frac{1}{t}\left(\int_O\varphi^*_{\exp tY}\circ
\phi^*_g\circ
A_g\,dg -\int_O\phi^*_g\circ A_g\,dg\right)\circ\psi(D_{\FF})\\
&=\lim_{t\to 0}\frac{1}{t}\left(\int_O\phi^*_{g}\circ
F^*_{\exp tY,g}\circ A_{\exp tY\cdot g}\,dg -\int_O\phi^*_g\circ A_g\,dg\right)\circ\psi(D_{\FF})\\
&=\int_O\phi^*_{g}\circ(V_Y\circ
A+L_YA)_g\,dg\circ\psi(D_{\FF})\;.
\end{align*}
So $\widehat{Y}\circ P\in\DD$.
\end{proof}

With the above notation, by the proof of
Proposition~\ref{p:smoothing} and~\eqref{e:fdnorm1}, it can be
easily seen that, for integers $k\le\ell$, there are some $C,C'>0$
and some natural $N$ such that
\begin{equation}\label{e:|P| k,ell}
\|P\|_{k,\ell}\leq C'\int
|(\id-\partial^2_{\xi})^N\hat{\psi}(\xi)|\,e^{C|\xi|}\,d\xi\;.
\end{equation}
Here, $C$ depends on $k$ and $\ell$, and $C'$ depends on $k$,
$\ell$ and $A$.

\subsection{A norm estimate}\label{ss:DD 0}

Let
$$
P=\int_O\phi_g^*\cdot f(g)\,\Lambda(g)\circ\psi(D_\FF)\in\DD\;,
$$
where $\phi$ and $\psi$ are like in Section~\ref{ss:DD}, and
$f\in\cinfc(O)$. In this case,~\eqref{e:|P| k,ell} is improved by
the following result, where $\D_G$ denotes the Laplacian of $G$.

\begin{prop}\label{p:|P| k,ell}
Let $K\subset O$ be a compact subset containing $\supp f$. For
naturals $k\le\ell$, there are some $C,C''>0$ and some natural
$N$, depending only on $K$, $k$ and $\ell$, such that
$$
\|P\|_{k,\ell} \leq C''\,\max_{g\in K}|(\id+\Delta_G)^Nf(g)|\,\int
|(\id-\partial^2_{\xi})^N\hat{\psi}(\xi)|\,e^{C|\xi|}\,d\xi\;.
$$
\end{prop}

\begin{proof}
Fix an orthonormal frame $Y_1,\dots,Y_q$ of $\mathfrak g$.
Consider any multi-index $J=(j_1,\dots,j_k)$ with
$j_1,\dots,j_k\in\{1,\dots,q\}$. We use the standard notation
$|J|=k$, and, with the notation of the proof of
Proposition~\ref{p:smoothing}, let:
\begin{itemize}

\item  $Y_J=Y_{j_1}\circ\cdots\circ Y_{j_k}$ (operating in $\cinf(G)$);

\item $\widehat{Y}_J=\widehat{Y_{j_1}}\circ\dots\circ\widehat{Y_{j_k}}$;

\item $V_J=V_{Y_{j_1}}\circ\cdots\circ V_{Y_{j_k}}$; and

\item $L_JA=L_{Y_{j_1}}\cdots L_{Y_{j_k}}A$ for any smooth family $A$ of
leafwise differential operators in $\Omega(\FF)$.

\end{itemize}
Consider the empty multi-index $\emptyset$ too, with
$|\emptyset|=0$, and define:
\begin{itemize}

\item  $Y_\emptyset=\id_{\cinf(G)}$;

\item $\widehat{Y}_\emptyset=\id_{\Omega(\FF)}$;

\item $V_{\emptyset,g}=\id_{\Omega(\FF)}$ for all $g\in O$,
defining a smooth family $V_\emptyset$; and

\item $L_{\emptyset}A=A$ for any smooth family $A$ of leafwise
differential operators in $\Omega(\FF)$.

\end{itemize}

Given any natural $N$, there is some $C_1>0$ such that
\begin{gather*}
\|\phi_g^*\|_k\le C_1\;,\qquad \|(L_JV_{J'})_g\|\le C_1\;,\\
\|(Y_Jf)(g)\|\le C_1\,\max_{g\in K}|(\id+\D_G)^Nf(g)|\;,\\
\|(\id+\phi_g^{*-1}\circ\Delta_\FF\circ
\phi^*_g)^N\circ\psi(\D_\FF)\|_k \le
C_1\,\|(\id+\Delta_\FF)^N\circ\psi(D_\FF)\|_k
\end{gather*}
for all $g\in K$ and all multi-indices $J$ and $J'$ with
$|J|,|J'|\le N$.

For any multi-index $J$, we have $$ \widehat{Y}_J\circ
P=\int_O\phi_g^*\circ A_{J,g}\,\Lambda(g)\circ \psi(D_\FF)\;, $$
where $A_J=\{A_{J,g}\ |\ g\in G\}$ is the smooth family of
leafwise differential operators inductively defined by setting
\begin{align*}
A_{\emptyset,g}&=\id_{\Omega(\FF)}\cdot f(g)\;,\\
A_{(j,J)}&=V_j\circ A_J+L_jA_J\;.
\end{align*}

By induction on $|J|$, we easily get that $A_J$ is a sum of smooth
families of leafwise differential operators of the form
$$
L_{J_1}V_{J'_1}\circ\dots\circ L_{J_\ell}V_{J'_\ell}\cdot
Y_{J''}f\;,
$$
where $J_1,J'_1,\dots,J_\ell,J'_\ell,J''$ are possibly empty
multi-indices satisfying
$$
|J_1|+|J'_1|+\dots+|J_\ell|+|J'_\ell|+|J''|=|J|\;.
$$
So there is some $C_2>0$ such that
$$
\|A_{J,g}\|_k\le C_2\,\max_{g\in K}|(\id+\Delta_G)^Nf(g)|
$$
for all $g\in K$ and every multi-index $J$ with $|J|\le N$. Hence
\begin{align*}
\|\widehat{Y}_J\circ P\|_k&\le\int_O\|\phi_g^*\|_k\| A_{J,g}\|_k\,dg\,\|\psi(D_{\FF})\|_k\\
&\le C_1C_2\,\max_{g\in
K}|(\id+\D_G)^Nf(g)|\,\int|\hat{\psi}(\xi)|\,e^{C|\xi|}\,d\xi
\end{align*}
for some $C>0$ by~\eqref{e:fdnorm}. On the other hand,
\begin{align*}
\|(\id+\D_\FF)^N\circ P\|_k&\le\int_O\|(\id+\phi_g^{*-1}\circ\Delta_\FF\circ\phi^*_g)^N\circ \psi(\D_\FF)\|_k\,|f(g)|\,\Lambda(g)\\
&\le C_1\,\int_O\|(\id+\D_\FF)^N\circ\psi(\D_\FF)\|_k\,|f(g)|\,\Lambda(g)\\
&\le C_1\,\max_{g\in K}|f(g)|\,\int
|(\id-\partial^2_{\xi})^N\hat{\psi}(\xi)|\,e^{C|\xi|}\,d\xi
\end{align*}
for some $C>0$ by~\eqref{e:fdnorm1}. Now, the result follows
because $-\sum_{j=1}^q\widehat{Y_j}^2$ is transversely elliptic,
and $\D_\FF$ is leafwise elliptic.
\end{proof}

\subsection{Parameter independence of the supertrace}

Choose an even function in $\AAA$, which can be written as
$x\mapsto\psi(x^2)$. Take also a \cinf\ local  representation
$\phi:M\times O\to M$ of $\Phi$ and some $f\in C^\infty_c(O)$.
Then consider the one parameter family of operators $P_t\in\DD$,
$t>0$, defined by
$$
P_t=\int_O\phi_g^*\cdot f(g)\,\Lambda(g)\circ\psi(t\D_\FF)^2\;.
$$

\begin{lem}\label{l:Tr sP t}
$\Tr^sP_t$ is independent of $t$.
\end{lem}

\begin{proof}
The proof is similar to the proof of the corresponding result in
the heat equation proof of the Lefschetz trace formula (see {\it e.g.\/}
\cite{Roe98}). We have
\begin{align*}
\frac{d}{dt}\,\Tr^sP_t&
=2\Tr^s\int_O\phi_g^*\cdot f(g)\,\Lambda(g)\circ\D_\FF\circ\psi'(t\D_\FF)\circ\psi(t\D_\FF)\\
&=2\Tr\int_O\phi_g^*\cdot f(g)\,\Lambda(g)
\circ d_\FF^-\circ\delta_\FF^+\circ\psi'(t\D_\FF^+)\circ\psi(t\D_\FF^+)\\
&\phantom{=\text{}}\text{}-2\Tr\int_O\phi_g^*\cdot
f(g)\,\Lambda(g)
\circ d_\FF^+\circ\delta_\FF^-\circ\psi'(t\D_\FF^-)\circ\psi(t\D_\FF^-)\\
&\phantom{=\text{}}\text{}+2\Tr\int_O\phi_g^*\cdot
f(g)\,\Lambda(g)
\circ\delta_\FF^-\circ d_\FF^+\circ\psi'(t\D_\FF^+)\circ\psi(t\D_\FF^+)\\
&\phantom{=\text{}}\text{}-2\Tr\int_O\phi_g^*\cdot
f(g)\,\Lambda(g)\circ \delta_\FF^+\circ
d_\FF^-\circ\psi'(t\D_\FF^-)\circ\psi(t\D_\FF^-)\;.
\end{align*}
On the other hand, since the function $x\mapsto\psi'(x^2)$ is in
$\AAA$, we have
\begin{multline*}
\Tr\int_O\phi_g^*\cdot f(g)\,\Lambda(g)\circ d_\FF^\mp\circ\delta_\FF^\pm\circ\psi'(t\D_\FF^\pm)\circ\psi(t\D_\FF^\pm)\\
\begin{aligned}
&= \Tr d_\FF^\mp\circ\int_O\phi_t^*\cdot f(g)\,\Lambda(g)
\circ\psi'(t\D_\FF^\pm)\circ\psi(t\D_\FF^\pm)\circ \delta_\FF^\pm\\
&=\Tr\psi(t\D_\FF^\pm)\circ\delta_\FF^\pm\circ d_\FF^\mp\circ
\int_O\phi_g^*\cdot f(g)\,\Lambda(g)\circ\psi'(t\D_\FF^\pm)\\
&=\Tr\int_O\phi_g^*\cdot f(g)\,\Lambda(g)
\circ\psi'(t\D_\FF^\pm)\circ\psi(t\D_\FF^\pm)\circ\delta_\FF^\pm\circ d_\FF^\mp\\
&=\Tr\int_O\phi_g^*\cdot f(g)\,\Lambda(g) \circ\delta_\FF^\pm\circ
d_\FF^\mp\circ\psi'(t\D_\FF^\pm)\circ\psi(t\D_\FF^\pm)\;,
\end{aligned}
\end{multline*}
where we have used the well known fact that, if $A$ is a trace
class operator and $B$ is bounded, then $AB$ and $BA$ are trace
class operators with the same trace. Therefore
$\frac{d}{dt}\,\Tr^sP_t=0$ as desired.
\end{proof}

\subsection{The global action on the leafwise complex}\label{ss:global action}

Let $\mathfrak{G}$ be the holonomy groupoid of $\FF$. Since the
leaves of Lie foliations have trivial holonomy groups, we have $$
\mathfrak{G}\equiv \{(x,y)\in M\times M \ |\  \text{$x$ and $y$
lie in the same leaf of $\FF$}\}\;. $$ This is a $\cinf$
submanifold of $M\times M$ which contains the diagonal $\D_M$. Let
$d_\FF$ be the distance function of the leaves of $\FF$. For each
$r>0$, the {\em $r$-penumbra\/} of $\D_M$ in $\mathfrak G$ is
defined by $$ \Pen_{\mathfrak G}(\D_M,r)=\{(x,y)\in\mathfrak G\ |\
d_\FF(x,y)< r\}\;. $$ Observe that a subset of $\mathfrak{G}$ has
compact closure if and only if it is contained in some penumbra of
$\D_M$. The product of two elements
$(x_1,y_1),(x_2,y_2)\in\mathfrak{G}$ is defined when $y_1=x_2$,
and it is equal to $(x_1,y_2)$. The space of units of
$\mathfrak{G}$ is $\D_M\equiv M$. The source and target
projections $s,r:\mathfrak{G}\to M$ are the restrictions of the
first and second factor projections $M\times M\to M$; thus $$
r^{-1}(x)=L_x\times\{x\}\;,\quad s^{-1}(x)=\{x\}\times L_x $$ for
each $x\in M$.

Let $S$ denote the \cinf\ vector bundle
$$
s^*\bigwedge T\FF^*\otimes r^*\bigwedge T\FF
$$
over $\mathfrak{G}$; thus
$$
S_{(x,y)}\equiv\bigwedge T_x\FF^*\otimes\bigwedge T_y\FF
\equiv\Hom(\bigwedge T_y\FF^*,\bigwedge T_x\FF^*)
$$
for each $(x,y)\in\mathfrak{G}$. Let $\omega_\FF$ be the volume form
of the leaves of $\FF$ (we assume that $\FF$ is oriented). Recall
that $\cinfc(S)$ is an algebra with the convolution product given by
$$
(k_1\cdot k_2)(x,y)=\int_{L_x}k_1(x,z)\circ
k_2(z,y)\,\omega_\FF(z)
$$
for $k_1,k_2\in\cinfc(S)$ and $(x,y)\in\mathfrak{G}$. Recall also
that the {\em global action\/} of $\cinfc(S)$ in $\Omega(\FF)$ is
defined by
$$
(k\cdot\alpha)(x)=\int_{L_x}k(x,y)\,\alpha(y)\,\omega_{\FF}(y)
$$
for $k\in\cinfc(S)$, $\alpha\in\Omega(\FF)$ and $x\in M$.

Consider the lift to $\widetilde{M}$ of the bundle-like metric of
$M$, and its restriction to the leaves of $\widetilde{\FF}$. Let
$U\Omega(\widetilde{\FF})\subset\Omega(\widetilde{\FF})$ be the
subcomplex of differential forms $\alpha$ whose covariant
derivatives $\nabla^r\alpha$ of arbitrary order $r$ are uniformly
bounded; this is a Frechet space with the metric induced by the
seminorms
$$
\|\|\alpha\|\|_r=\sup\{\nabla^r\alpha(\tilde x)\ |\ \tilde
x\in\widetilde{M}\}\;.
$$
Observe that $\pi^*(\Omega(\FF))\subset U\Omega(\widetilde{\FF})$.

The holonomy groupoid $\widetilde{\mathfrak{G}}$ of
$\widetilde{\FF}$ satisfies the same properties as $\mathfrak{G}$,
except that, in $\widetilde{\mathfrak{G}}$, the penumbras of the
diagonal $\D_{\widetilde{M}}$ have compact closure if and only
$\widetilde{M}$ is compact.

The map $\pi\times\pi:\widetilde{M}\times \widetilde{M}\to M\times
M$ restricts to a covering map
$\widetilde{\mathfrak{G}}\to\mathfrak{G}$, whose group of deck
transformations is isomorphic to $\Aut(\pi)$: for each
$\sigma\in\Aut(\pi)$, the corresponding element in
$\Aut(\widetilde{\mathfrak{G}}\to\mathfrak{G})$ is the restriction
$\sigma\times\sigma:\widetilde{\mathfrak{G}}\to\widetilde{\mathfrak{G}}$.

Let $\widetilde{S}$ denote the \cinf\ vector bundle
$$
\tilde s^*\bigwedge T\widetilde{\FF}^*\otimes\tilde r^*\bigwedge T\widetilde{\FF}
$$
over $\widetilde{\mathfrak{G}}$, and let
$C^\infty_\D(\widetilde{S})\subset\cinf(\widetilde{S})$ denote the
subspace of sections supported in some penumbra of
$\D_{\widetilde{M}}$. As above, this set becomes an algebra with
the convolution product, and there is a {\em global action\/} of
$C^\infty_\D(\widetilde{S})$ in $U\Omega(\widetilde{\FF})$.

Any $k\in\cinf(S)$ lifts via $\pi\times\pi$ to a section $\tilde
k\in\cinf(\widetilde{S})$. Since $\pi$ restricts to
diffeomorphisms of the leaves of $\widetilde{\FF}$ to the leaves
of $\FF$, it follows that $\tilde k\in C^\infty_\D(\widetilde{S})$
if $k\in\cinfc(S)$.

Take any $\psi\in\AAA$. For each leaf $L$ of $\FF$, denoting by
$\D_L$ the Laplacian of $L$, the spectral theorem defines a
smoothing operator $\psi(\D_L)$ in $\bfOmega(L)$, and the family
$$
\{\psi(\D_L)\ |\ \text{$L$ is a leaf of $\FF$}\}
$$
is also denoted by $\psi(\D_\FF)$. By
\cite[Proposition~2.10]{Roe88I}, the Schwartz kernels $k_L$ of
the operators $\psi(\D_L)$
can be combined to define a section $k\in\cinf(S)$, called the {\em
leafwise smoothing kernel\/} or {\em leafwise Schwartz kernel\/} of
$\psi(\D_\FF)$.

Suppose that the Fourier transform $\hat\psi$ of $\psi$ is
supported in $[-R,R]$ for some $R>0$. Then, according to the proof
of Assertion~1 in \cite[page~461]{Roe87}, $k$ is supported in the
$R$-penumbra of $\D_M$, and thus $k\in\cinfc(S)$. Moreover the
operator $\psi(D_\FF)$ in $\Omega(\FF)$, defined by the spectral
theorem, equals the operator given by the global action of $k$.

Consider also the lift $\tilde k\in C^\infty_\D(\widetilde S)$,
whose global action in $U\Omega(\widetilde{\FF})$ defines an
operator denoted by $\psi(D_{\widetilde{\FF}})$. It is clear that
the diagram
\begin{equation}\label{e:psi(D widetilde FF)}
\begin{CD}
U\Omega(\widetilde{\FF}) & @>{\psi(D_{\widetilde{\FF}})}>> &U\Omega(\widetilde{\FF})\\
@A{\pi^*}AA & & @AA{\pi^*}A\\
\Omega(\FF) &@>{\psi(D_{{\FF}})}>> &\Omega(\FF)
\end{CD}
\end{equation}
commutes.

Any function $\psi\in\AAA$ with compactly supported Fourier
transform can be modified as follows to achieve the condition of
being supported in $[-R,R]$. For each $t>0$, let $\psi_t\in\AAA$
be the function defined by $\psi_t(x)=\psi(tx)$.

\begin{lem}\label{l:psi t}
If $\hat\psi$ is compactly supported for some $\psi\in\AAA$, then
$\widehat{\psi_t}$ is supported in $[-R,R]$ for $t$ small enough.
\end{lem}

\begin{proof}
This holds because
$\widehat{\psi_t}(\xi)=\frac{1}{t}\hat\psi(\frac{\xi}{t})$.
\end{proof}

\subsection{Schwartz kernels}\label{ss:kernels}

Let $\phi$, $f$, $\psi$ and $P$ be like in Section~\ref{ss:DD 0}
such that $\hat\psi$ is compactly supported. Take some $R>0$ so that
$\supp\hat\psi\subset[-R,R]$. Let $k\in\cinfc(S)$ be the leafwise kernel of $\psi(D_\FF)$, and let
$\tilde{k}\in C^\infty_\D(\widetilde S)$ be the lift of $k$, whose
action in $\Omega(\widetilde\FF)$ defines the operator
$\psi(D_{\widetilde\FF})$ (Section~\ref {ss:global action}).

Let $\tilde\phi:\widetilde{M}\times O\to\widetilde{M}$ be a \cinf\
lift of $\phi$. Define $\widetilde{P}:U\Omega(\widetilde{\FF})\to
U\Omega(\widetilde{\FF})$ by $$ \widetilde{P}=
\int_O\tilde{\phi}^*_g\cdot
f(g)\,\Lambda(g)\circ\psi(D_{\widetilde{\FF}}). $$  The
commutativity of the diagram $$
\begin{CD}
U\Omega(\widetilde{\FF}) & @>{\widetilde{P}}>> &U\Omega(\widetilde{\FF})\\
@A{\pi^*}AA & & @AA{\pi^*}A\\
\Omega(\FF) &@>P>> &\Omega(\FF)
\end{CD}
$$
follows from the commutativity of \eqref{e:psi(D widetilde FF)}.

Let $\omega_{\widetilde\FF}$ be the volume form of the leaves of
$\widetilde\FF$, which can be also considered as a differential form
on $M$ that vanishes when some vector is orthogonal to the leaves.
Thus the volume form of $\widetilde M$ is $\omega_{\widetilde
M}=D^*\Lambda\wedge\omega_{\widetilde\FF}$ with the right choice of
orientations. For $\tilde{x}\in\widetilde{M}$ and $\alpha\in
U\Omega(\widetilde{\FF})$, we have
\begin{align*}
(\widetilde{P}\alpha)(\tilde{x})
&=(\int_O\tilde{\phi}^*_g \cdot f(g)\,\Lambda(g)\circ \psi(D_{\widetilde{\FF}})\alpha)(\tilde{x})\\
&=\int_O\tilde{\phi}^*_g((\psi(D_{\widetilde{\FF}})\alpha)(\tilde\phi_g(\tilde{x})) \cdot f(g)\,\Lambda(g)\\
&=\int_O\int_{\widetilde{L}_{\tilde{x}}}\tilde{\phi}^*_g\circ
\tilde{k}(\tilde\phi_g(\tilde{x}),\tilde{y})(\alpha(\tilde{y}))\,\omega_{\widetilde{\FF}}(\tilde{y}) \cdot f(g)\,\Lambda(g)\\
&=\int_{\phi(\widetilde{L}_{\tilde x}\times
O)}\tilde{\phi}^*_g\circ
\tilde{k}(\tilde\phi_g(\tilde{x}),\tilde{y})(\alpha(\tilde{y}))
\cdot f(g)\,\omega_{\widetilde M}(\tilde y)
\end{align*}
by Corollary~\ref{c:tilde varphi}, where $g\in O$ is determined by
the condition $\tilde y\in\tilde\phi_g(\widetilde L_{\tilde x})$,
which means $g=D(\tilde x)^{-1}D(\tilde y)$ by Lemma~\ref{l:tilde
varphi}. So we can say that $\widetilde{P}$ is  given by the
Schwartz kernel $\tilde p$ defined by
\begin{equation}\label{e:tilde p}
\tilde{p}(\tilde{x},\tilde{y})=
\begin{cases}
\tilde{\phi}^*_g\circ
\tilde{k}(\tilde\phi_g(\tilde{x}),\tilde{y})\cdot f(g)
& \text{if $\tilde{y}\in \tilde{\phi}(\widetilde{L}_{\tilde{x}}\times O)$}\\
0 & \text{otherwise}
\end{cases}
\end{equation}
for $g\in O$ as above. It follows that
\begin{equation}\label{e:p}
p(x,y)=\sum_{\sigma\in\operatorname{Aut}(\pi)}\tilde{p}(\tilde{x},\sigma(\tilde{y}))\;,
\end{equation}
where $\tilde{x}\in\pi^{-1}(x)$, $\tilde{y}\in\pi^{-1}(y)$, and we
use identifications $T_{\tilde{x}}\widetilde{\FF}\equiv T_x\FF$
and $T_{\sigma(\tilde{y})}\widetilde{\FF}\equiv T_y\FF$ given by
$\pi_*$.

For each $x\in M$, $\tilde x\in\widetilde M$ and $r>0$, let
$B_\FF(x,r)$ and $B_{\widetilde\FF}(\tilde x,r)$ be the $r$-balls
of centers $x$ and $\tilde x$ in $L_x$ and $\widetilde L_{\tilde
x}$, respectively. Let $O_1$ be an open subset of $G$ whose
closure is compact and contained in $O$. By the compactness of
$M\times\overline O_1$, there is some $R_1>0$ such that
\begin{equation}\label{e:R 1}
B_\FF(\phi_g(x),R)\subset\phi_g(B_\FF(x,R_1))
\end{equation}
for all $x\in M$ and all $g\in O_1$. So
\begin{equation}\label{e:R 1, tilde}
B_{\widetilde\FF}(\tilde\phi_g(\tilde
x),R)\subset\tilde\phi_g(B_{\widetilde{\FF}}(\tilde{x},R_1))
\end{equation}
for all $\tilde x\in\widetilde M$ and all $g\in O_1$ because $\pi$
restricts to isometries of the leaves of $\widetilde\FF$ to the
leaves of $\FF$.

\begin{lem}\label{l:injective}
Each $g\in O$ has a neighborhood $O_1$ as above such that
$$
\pi:\tilde\phi(\overline{B_{\widetilde{\FF}}(\tilde{x},R_1)}\times
O_1)\to M
$$
is injective for any $\tilde{x}\in\widetilde{M}$.
\end{lem}

\begin{proof}
Since $M$ is compact, there exists a compact subset
$K\subset\widetilde{M}$ with $\pi(K)=M$. Notice that, if the
statement holds for some $\tilde{x}\in \widetilde{M}$, then it also
holds for all points in the $\Aut(\pi)$-orbit of $\tilde x$. So, if
the statement fails, there exist sequences
$\tilde{x}_i,\tilde{y}_i\in \widetilde{M}$ and $\sigma_i\in
\Aut(\pi)$ such that $\tilde{x}_i\in K$,
$\sigma_i\neq\id_{\widetilde{M}}$, and $$ d_{\widetilde M
}(\{\tilde{y}_i,\sigma_i(\tilde{y}_i)\},\tilde\phi_g(B_{\widetilde{\FF}}(\tilde{x}_i,R_1)))\to0
$$
as $i\to\infty$; observe that $D(\tilde x_i)^{-1}\,D(\tilde y_i)\to
g$ by Lemma~\ref{l:tilde varphi}. Since $K$ is compact, we can
assume that there exists $\lim_i
\tilde{x}_i=\tilde{x}\in\widetilde{M}$, where $d_{\widetilde M}$
denotes the distance function of $\widetilde M$. Hence $\tilde{y}_i$
and $\sigma_i(\tilde{y}_i)$ approach
$\tilde\phi_g(B_{\widetilde{\FF}}(\tilde{x},R_1))$. Since
$\tilde\phi_g(B_{\widetilde{\FF}}(\tilde{x},R_1))$ has compact
closure, it follows that $\tilde{y}_i$ and $\sigma_i(\tilde{y}_i)$
lie in some compact neighborhood $Q$ of
$\tilde\phi_g(B_{\widetilde{\FF}}(\tilde{x},R_1))$ for infinitely
many indices $i$, yielding $\sigma_i(Q)\bigcap Q\not=\emptyset$. So
there is some $\sigma\in\operatorname{Aut}(\pi)$ such that
$\sigma_i=\sigma$ for infinitely many indices $i$. In particular,
$\sigma\neq\id_{\widetilde{M}}$.

On the other hand, since  $\tilde{y}_i$ and
$\sigma_i(\tilde{y}_i)$ approach
$\tilde\phi_g(B_{\widetilde{\FF}}(\tilde x,R_1))$, which has
compact closure, we can assume that there exist $\lim_i
\tilde{y}_i=\tilde{y}$ and $\lim_i
\sigma_i(\tilde{y}_i)=\sigma(\tilde y)$ in
$\overline{\tilde\phi_g(B_{\widetilde{\FF}}(\tilde x,R_1))}$,
which is contained in the leaf $\tilde\phi_g(\widetilde{L}_{\tilde
x})$ (a fiber of $D$). So
$$
D(\tilde{y})=D(\sigma(\tilde{y}))=h(\sigma)\cdot D(\tilde{y})\;,
$$
yielding $h(\sigma)=e$, and thus $\sigma=\id_{\widetilde{M}}$
because $h$ is injective. This contradiction concludes the proof.
\end{proof}

From now on, assume that $\phi$ satisfies~\eqref{e:R 1} and the
property of the statement of Lemma~\ref{l:injective} with some
fixed open subset $O_1\subset O$ which contains the support of
$f$.

\begin{cor}\label{c:injective}
The map $\pi$ is injective on the support of
$\tilde{p}(\tilde{x},\cdot)$ for any $\tilde{x}\in\widetilde{M}$.
\end{cor}

\begin{proof}
By~\eqref{e:tilde p},~\eqref{e:R 1, tilde} and since $\tilde k$ is
supported in the $R$-penumbra of $\Delta_{\widetilde M}$, we get
$$
\supp(\tilde
p(\tilde{x},\cdot))\subset\tilde\phi(B_{\widetilde{\FF}}(\tilde{x},R_1)\times
O_1)
$$
for any $\tilde{x}\in\widetilde{M}$, and the result follows from
Lemma~\ref{l:injective}.
\end{proof}

\begin{cor}\label{c:p}
We have
$$
p(x,y)=
\begin{cases}
\phi^*_g\circ k(\phi_g(x),y)\cdot f(g)
& \text{if $y\in\phi(L_x\times O_1)$}\\
0 & \text{otherwise,}
\end{cases}
$$
where $g\in O_1$ is determined by the condition
$y\in\phi_g(B_\FF(x,R_1))$.
\end{cor}

\begin{proof}
This is a consequence of~\eqref{e:tilde p},~\eqref{e:p},
Corollary~\ref{c:injective} and Lemma~\ref{l:injective}.
\end{proof}

\begin{cor}\label{c:p(x,x)}
If $e\in O_1$ and $\phi_e=\id_M$, then
$$
p(x,x)=k(x,x)\cdot f(e)\;.
$$
\end{cor}

\begin{proof}
Since $\phi_e=\id_M$, the result follows from Corollary~\ref{c:p}
and the following assertion.

\begin{claim}\label{cl:g=e}
For all $g\in O_1$ and $x\in M$, if $x\in\phi_g(B_\FF(x,R_1))$,
then $g=e$.
\end{claim}

By Lemma~\ref{l:injective},
$$
\pi:\tilde\phi(B_{\widetilde{\FF}}(\tilde x,R_1)\times
O_1)\to\phi(B_\FF(x,R_1)\times O_1)
$$
is a diffeomorphism. On the other hand,
$$
\tilde\phi:\widetilde{L}_{\tilde x}\times
O_1\to\tilde\phi(\widetilde{L}_{\tilde x}\times O_1)
$$
is a diffeomorphism as well by Corollary~\ref{c:tilde varphi}. It
follows that
$$
\phi:B_\FF(x,R_1)\times O_1\to\phi(B_\FF(x,R_1)\times O_1)
$$
is also a diffeomorphism, which implies Claim~\ref{cl:g=e} because
$\phi_e(x)=x$.
\end{proof}

\begin{lem}\label{l:g 2=a -1g 1a}
For $i\in\{1,2\}$, suppose that $x_i\in\phi_{g_i}(B_\FF(x_i,R_1))$
for some $(x_i,g_i)\in M\times O_1$. If $x_2$ is close enough to
$x_1$, then there is some $a\in G$ such that
$x_2\in\Phi_a(L_{x_1})$ and $g_2=a^{-1}g_1a$.
\end{lem}

\begin{proof}
We have
$$
\Phi_a(L_{x_1})=\Phi_a\circ\Phi_{g_1}(L_{x_1})=\Phi_{a^{-1}g_1a}\circ\Phi_a(L_{x_1})
$$
for all $a\in G$. Therefore, if $x_2$ is close enough to $x_1$,
there is some $a\in G$ such that $a^{-1}g_1a\in O_1$ and
$$
x_2\in\Phi_a(L_{x_1})\cap\phi_{a^{-1}g_1a}(B_\FF(x_2,R_1))\;.
$$
Then the result follows because the condition
$x_2\in\phi_{g_2}(B_\FF(x_2,R_1))$ determines $g_2$ in $O_1$ by
Lemma~\ref{l:injective}.
\end{proof}

\section{Lefschetz distribution}\label{s:Lefschets dis}

Let $\phi:M\times O\to M$ be a \cinf\ local representation of the
structural transverse action $\Phi$ on some open subset $O\subset
G$. For any $f\in \cinfc(O)$ and $t>0$, let $P_f$ and $Q_{t,f}$ be
the operators in $\Omega(\FF)$ defined by
\begin{align*}
P_f&=\int_O\phi^*_g\cdot f(g)\,\Lambda(g)\circ\Pi\;,\\
Q_{t,f}&=\int_O\phi^*_g\cdot f(g)\,\Lambda(g)\circ
e^{-t\Delta_\FF}\;.
\end{align*}
The operator $Q_{t,f}$ is in the class $\DD$, and thus it is
smoothing by Proposition~\ref{p:smoothing}.

\begin{prop}\label{p:P f is smoothing}
$P_f$ is a smoothing operator.
\end{prop}

\begin{proof}
By \cite{AlvKordy:heat}, $\Pi$ defines a bounded operator in each
Sobolev space $W^k\Omega^i(\FF)$. Hence,  $P_f=Q_{t,f}\circ \Pi$
is smoothing because so is $Q_{t,f}$.
\end{proof}

By Proposition~\ref{p:P f is smoothing}, $P_f$ is a trace class
operator in the space $\Omega(\FF)$, and thus so is $P_f^{(i)}$.

\begin{prop}\label{p:Tr P f (i)}
The functional $f\mapsto \Tr P_f^{(i)}$ is a distribution on $O$.
\end{prop}

\begin{proof}
Since $\bfPi$ is a projection in $\bfOmega(\FF)$ and
$P_f=Q_{t,f}\circ\Pi$, we have
$$
\|P_f^{(i)}\|_{0,k}\leq\|Q_{t,f}^{(i)}\|_{0,k}\;,
$$
and the result follows by~\eqref{e:Tr} and Proposition~\ref{p:|P|
k,ell}.
\end{proof}

Proposition~\ref{p:trace class} is given by Propositions~\ref{p:P
f is smoothing} and~\ref{p:Tr P f (i)}.

Because the endomorphism $\Phi^*_g$ of $\overline H(\FF)$
corresponds to the operator $\Pi\circ\phi_g^*$ in $\HH(\FF)$ by the
leafwise Hodge isomorphism, the composite $\Pi\circ P_f$ is
independent of the choice of $\phi$. Moreover $\Tr
P_f^{(i)}=\Tr(\Pi\circ P_f^{(i)})$. Hence the distributions given by
Proposition~\ref{p:Tr P f (i)} can be combined to form a global
distribution $\Tr_{\text{\rm dis}}^i(\FF)$ on $G$; in this notation,
$\FF$ refers to the foliation endowed with the given transverse Lie
structure, which indeed is determined by the foliation when the
leaves are dense. Each $\Tr_{\text{\rm dis}}^i(\FF)$ is called a
{\em distributional trace\/} of $\FF$, and define the {\em Lefschetz
distribution\/} of $\FF$ by the formula
$$
L_{\text{\rm dis}}(\FF)=\sum_i (-1)^i\,\Tr_{\text{\rm
dis}}^i(\FF)\;.
$$

\begin{lem}\label{l:converge}
For any $f\in\cinfc(O)$, $\Tr Q^{(i)}_{t,f}\to\Tr P^{(i)}_f$ as
$t\to\infty$.
\end{lem}

\begin{proof}
Since $Q_{1,f}$ is smoothing, it defines a bounded operator
$W^{-1}\Omega^i(\FF)\to W^k\Omega^i(\FF)$ for any $k$. By
\cite{AlvKordy:heat}, $e^{-(t-1)\Delta_\FF}-\Pi$ is bounded in
$W^{-1}\Omega^i(\FF)$ for $t>1$ and converges strongly to $0$ as
$t\rightarrow\infty$. From the compactness of the canonical
embedding $\bfOmega^i(\FF)\to W^{-1}\Omega^i(\FF)$, it follows
that $e^{-(t-1)\Delta_\FF}-\Pi$ converges uniformly to $0$ as
$t\rightarrow \infty$ as an operator $\bfOmega^i(\FF)\to
W^{-1}\Omega^i(\FF)$. Therefore $\|Q_{t,f}-P_f\|_{0,k}\to0$ as
$t\rightarrow \infty$ for any $k$ because
$$
Q_{t,f}-P_f=Q_{1,f}\circ(e^{-(t-1)\Delta_\FF}-\Pi)\;.
$$
Then the result follows from~\eqref{e:Tr}.
\end{proof}

\begin{cor}\label{c:Tr sQ t,f}
$\Tr^sQ_{t,f}=\Tr^s P_f$ for all $t$.
\end{cor}

\begin{proof}
This follows from Lemmas~\ref{l:Tr sP t} and~\ref{l:converge}.
\end{proof}

\section{The distributional Gauss-Bonett theorem}\label{s:Gauss-Bonett}

The holonomy pseudogroup of $\FF$ is represented by the
pseudogroup on $G$ generated by the left translations given by
elements of $\Gamma$. Thus $\Lambda$ can be considered as a
holonomy invariant transverse measure of $\FF$. To be more
precise, take a $(G,\Gamma)$-valued foliated cocycle $\{U_i,f_i\}$
defining the given transverse Lie structure (Section~\ref{s:Lie}).
The differential forms $f_i^*\Lambda$ can be combined to get the
transverse volume form $\omega_\Lambda$ of $\FF$. We can also
describe $\omega_\Lambda$ by the condition
$D^*\Lambda=\pi^*\omega_\Lambda$. The restriction of
$\omega_\Lambda$ to smooth local transversals is the precise
interpretation of $\Lambda$ as a holonomy invariant measure on
local transversals.

By non-commutative integration theory \cite{Co79}, the holonomy
invariant transverse measure $\Lambda$ defines a trace $\Tr_\Lambda$
on the twisted foliation von Neumann algebra $W^*(M,\FF,\bigwedge
T\FF^*)$. Consider also the corresponding supertrace
$\Tr_\Lambda^s$, equal to $\pm\Tr_\Lambda$, depending on whether the
even-odd bigrading is preserved or interchanged.

With the notation of Section~\ref{ss:global action}, we have
$\cinfc(S)\subset W^*(M,\FF,\bigwedge T\FF^*)$; here, each
$k\in\cinfc(S)$ is identified to the family of operators on the
leaves whose Schwartz kernels are the restrictions of $k$, and
moreover
$$
\Tr_{\Lambda}(k)=\int_M \Tr k(x,x)\,\omega_M(x)\;,\quad
\Tr_{\Lambda}^s(k)=\int_M \Tr^s k(x,x)\,\omega_M(x)\;.
$$

For each leaf $L$, let $\bfOmega(L)$ denote the Hilbert space of
$L^2$ differential forms on $L$, let
$\HH(L)\subset\bfOmega(L)$ be the subspace of harmonic $L^2$ forms,
and let $\Pi_L$ be the orthogonal projection $\bfOmega(L)\to\HH(L)$.
The family
$$
\Pi_\FF=\{\Pi_L \ |\  \text{$L$ is a leaf of $\FF$}\}
$$
defines a projection in $W^*(M,\FF,\bigwedge T\FF^*)$. The notation
$\Pi_L^{(i)}$ and $\Pi_\FF^{(i)}$ is used when we are only
considering differential forms of degree $i$. For each leaf $L$, let
$S_L=S|_{L\times L}$, and let $k_L,k_L^{(i)}\in\cinf(S_L)$ denote
the Schwartz kernels of $\Pi_L$ and $\Pi_L^{(i)}$. These sections
can be combined to define measurable sections $k$ and $k^{(i)}$ of
$S$, called the {\em leafwise Schwartz kernels\/} of $\Pi_\FF$ and
$\Pi_\FF^{(i)}$. Since $k$ and $k^{(i)}$ are \cinf\ along the fibers
of the source and target projections, their restrictions to the
diagonal $\D_M$ are measurable, and we have
$$
\Tr_{\Lambda}(\Pi_\FF^{(i)})=\int_M \Tr
k^{(i)}(x,x)\,\omega_M(x)\;,\quad \Tr_{\Lambda}^s(\Pi_\FF)=\int_M \Tr^s
k(x,x)\,\omega_M(x)\;.
$$
According to~\cite{Co79}, the $i$th {\em $\Lambda$-Betti number\/}
is defined by
$$
\beta^i_{\Lambda}(\FF)=\Tr_{\Lambda}(\Pi_\FF^{(i)})\;,
$$
and the {\em $\Lambda$-Euler characteristic\/} is given by the
formula
$$
\chi_{\Lambda}(\FF)=\Tr_\Lambda^s(\Pi_\FF)=\sum_i(-1)^i\,\beta^i_{\Lambda}(\FF)\;.
$$

\begin{thm}\label{t:L dis(FF)=chi Lambda(FF) cdot delta e}
$L_{\text{\rm dis}}(\FF)=\chi_{\Lambda}(\FF)\cdot\delta_e$ in some
neighborhood of $e$.
\end{thm}

Like in \cite[p.~463]{Roe87}, choose a sequence of smooth even
functions on $\R$, written as $x\mapsto \psi_m(x^2)$ with
$\psi_m(0)=1$, whose Fourier transforms are compactly supported and
which tend to the function $x\mapsto e^{-x^2/2}$ in the Schwartz
space $\SSS(\R)$. Let $k_{m,t}$ be the leafwise Schwartz kernel of
$\psi_m(t\D_\FF)^2$, which is in $\cinfc(S)$ according to
\cite{Roe87}. In \cite[p.~463]{Roe87}, it is proved that
\begin{equation}\label{e:Tr s Lambda psi m(t D FF) 2}
\Tr^s_{\Lambda}\psi_m(t\D_\FF)^2=\chi_{\Lambda}(\FF)\;.
\end{equation}

Let $\phi:M\times O\to M$ be any \cinf\ local representation of
$\Phi$ on some neighborhood $O$ of
$e$ such that $\phi_e=\id_M$, whose existence is given by
Lemma~\ref{l:phi e=id M}. Given $R>0$, take $R_1>0$  and some open
subset $O_1$ of $O$ containing $e$ such that~\eqref{e:R 1} and
Lemma~\ref{l:injective} are satisfied.

For every $f\in\cinfc(O)$ supported in $O_1$, let
$$
Q_{m,t,f}=\int_O\phi_g^*\cdot
f(g)\,\Lambda(g)\circ\psi_m(t\Delta_\FF)^2\in\DD\;.
$$

\begin{lem}\label{l:Tr sQ m,t,f}
$\Tr^sQ_{m,t,f}=\chi_\Lambda(\FF)\cdot f(e)$.
\end{lem}

\begin{proof}
By Lemma~\ref{l:psi t}, we can apply Corollary~\ref{c:p(x,x)} to
$Q_{m,t,f}$ when $t$ is small enough, obtaining
\begin{align*}
\Tr^s Q_{m,t,f} & = \int_M\Tr^sk_{m,t}(x,x)\cdot f(e)\,\omega_M(x)\\
                        & = \Tr^s_{\Lambda}\psi_m(t\Delta_\FF)^2\cdot f(e)\;.
\end{align*}
Then the result follows by~\eqref{e:Tr s Lambda psi m(t D FF) 2}.
\end{proof}

Consider the operators $Q_{t,f}$ and $P_f$ of
Section~\ref{s:Lefschets dis}.

\begin{lem}\label{l:lim m to infty Tr s Q m,t,f}
We have
$$
\lim_{m\to\infty}\Tr^s Q_{m,t,f}=\Tr^s Q_{t,f}
$$
for each $t$.
\end{lem}

\begin{proof}
Since the function $x\mapsto\psi_m(tx^2)-e^{-\frac{t}{2}x^2}$ tends to zero in $\AAA$ as $m\to\infty$, we get
$$
\lim_{m\to\infty}\|Q_{m,t,f}-Q_{t,f}\|_{0,k}=0
$$
for all $k$ by~\eqref{e:|P| k,ell} (or Lemma~\ref{p:|P| k,ell}),
and the result follows from~\eqref{e:Tr}
\end{proof}

Theorem~\ref{t:L dis(FF)=chi Lambda(FF) cdot delta e} follows from
Lemmas~\ref{l:Tr sQ m,t,f} and~\ref{l:lim m to infty Tr s Q
m,t,f}, and Corollary~\ref{c:Tr sQ t,f}.

\section{The distributional Lefschetz trace formula}\label{s:Lefschetz}

Let $\FF'$ be the foliation of $M\times G$ whose leaves are the
sets $L\times\{g\}$ for leaves $L$ of $\FF$ and points $g\in G$.
Lemma~\ref{l:g 2=a -1g 1a} suggests the following definition: for
each $x\in M$ and $g\in G$, let
\[
M'_{(x,g)}=\bigcup_{a\in G}(\Phi_a(L_x)\times\{a^{-1}ga\})\;.
\]
Observe that $M'_{(x,e)}=M\times\{e\}$. Moreover
$M'_{(x_1,g_1)}=M'_{(x_2,g_2)}$ if and only if $(x_2,g_2)\in
M'_{(x_1,g_1)}$; thus these sets form a partition of $M\times G$.

\begin{prop}\label{p:GG}
The sets $M'_{(x,g)}$ are the leaves of a \cinf\ foliation $\GG$
on $M\times G$.
\end{prop}

\begin{proof}
Consider the canonical identity $T_{(x,g)}(M\times G)\equiv
T_xM\oplus T_gG$ for each $(x,g)\in M\times G$, and let
$\Ad:G\to\Aut(\mathfrak{g})$ denote the adjoint representation of
$G$. With the notation of Section~\ref{s:structural action},
consider the \cinf\ vector subbundles $\VV,\WW\subset T(M\times
G)$ given by
\begin{align*}
\VV_{(x,g)}&=\{(X^\nu(x),(X-\Ad_{g^{-1}}(X))(g)) \ |\  X\in\mathfrak{g}\}\;,\\
\mathcal{W}_{(x,g)}&=\VV_{(x,g)}+T_{(x,g)}\FF'\;.
\end{align*}
The distribution defined by $\VV$ is not completely integrable.
Nevertheless, since $[X^\nu,Y^\nu]-[X,Y]^\nu\in\mathfrak{X}(\FF')$
for all $X,Y\in\mathfrak{g}$, it follows that the distribution
defined by $\WW$ is completely integrable. Thus there is a \cinf\
foliation $\GG$ on $M\times G$ so that $T\GG=\WW$. It is easy to
check that the leaves of $\GG$ are the sets $M'_{(x,g)}$.
\end{proof}

Let $\pr_1$ and $\pr_2$ denote the first and second factor
projections of $M\times G$ onto $M$ and $G$, respectively.

\begin{prop}\label{p:M'}
For each leaf $M'$ of $\GG$, we have the following:
\begin{itemize}

\item[(i)] the restriction $\pr_1:M'\to M$ is a covering map; and

\item[(ii)] $\pr_2$ restricts to a fiber bundle map of $M'$ to
some orbit of the adjoint action of $G$ on itself.

\end{itemize}
\end{prop}

\begin{proof}
For any $x\in M$, there is some open neighborhood $P$ of $x$ in
$L_x$, and some local representation $\varphi:M\times O\to M$ of
$\Phi$ on some open neighborhood $O$ of $e$ such that $\varphi$
restricts to a diffeomorphism of $P\times O$ onto some neighborhood
$U$ of $x$.  For any $g\in G$ such that $(x,g)\in M'$, the set
$$
\widetilde{U}_g=\{(\varphi_a(y),a^{-1}ga)\ |\ y\in P,\ a\in O\}
$$
is an open neighborhood of $(x,g)$ in $M'$, and the restriction
$\pr_1:\widetilde{U}_g\to U$ is a diffeomorphism. Therefore
property~(i) follows.

It is clear that $\pr_2(M')$ is an orbit of the adjoint action of
$G$ on itself, and that $\pr_2:M'\to\pr_2(M')$ is a \cinf\
submersion; thus its fibers are \cinf\ submanifolds. If $(x,g)\in
M'$, it can be easily seen that
$$
\pr_2^{-1}(g)\cap
M'=\{(\phi_a(y),g)\ |\ y\in L_x,\ a\in G_g,\ \phi_g\in\Phi_g\}\;,
$$
where $G_g$ is the centralizer of $g$ in $G$. For $\varphi:M\times
O\to M$ as above, the set $O'=\{b^{-1}gb\ |\ b\in O\}$ is an open
neighborhood of $g$ in $\pr_2(M')$. Let $$ F:O'\times
(\pr_2^{-1}(g)\cap M')\to\pr_2^{-1}(O')\cap M' $$ be the map defined
by $$
F(b^{-1}gb;\varphi_a(y),g)=(\varphi_{b^{-1}ab}\circ\varphi_b(y),b^{-1}gb)
$$ for $y\in L_x$, $a\in G_g$ and $b\in O'$. It is easy to see
that $F$ is a \cinf\ diffeomorphism, which shows property~(ii).
\end{proof}

Observe that $\FF'$ is a subfoliation of $\GG$, and, for each leaf
$M'$ of $\GG$, the restriction $\FF'|_{M'}$ is equal to the lift of $\FF$ by $\pr_1:M'\to M$.

Let $\phi:M\times O\to M$ be any \cinf\ local representation of
$\Phi$. Given $R>0$, take $R_1>0$  and some open subset $O_1$ of
$O$ containing $e$ such that~\eqref{e:R 1} and
Lemma~\ref{l:injective} are satisfied. Let
$$
\SSS=\{(x,g)\in M\times O_1\ |\ x\in\phi_g(B_\FF(x,R_1))\}\;.
$$

\begin{prop}\label{p:SSS}
We have:
\begin{itemize}

\item[(i)] $\SSS$ is contained in a finite union of leaves of $\GG$; and

\item[(ii)] the restriction $\pr_1:\SSS\to M$ is injective.

\end{itemize}
\end{prop}

\begin{proof}
Property~(i) is a consequence of Lemma~\ref{l:g 2=a -1g 1a} and
the compactness of $M$. Property~(ii) follows from
Lemma~\ref{l:injective}.
\end{proof}

Let $\phi':M\times O\to M\times O$ be the \cinf\ diffeomorphism
defined by $\phi'(x,g)=(\phi(x,g),g)$. Observe that $\phi'$ is a
foliated map $\FF'|_{M\times O}\to\FF'|_{M\times O}$.

\begin{prop}\label{p:phi'}
Let $M'$ be a leaf of $\GG$. If $\phi'$ preserves some leaf of
$\FF'|_{M'\cap(M\times O)}$, then it preserves every leaf of
$\FF'|_{M'\cap(M\times O)}$.
\end{prop}

\begin{proof}
Take some point $(x,g)$ in a leaf $L'$ of $\FF'|_{M'\cap(M\times
O)}$; thus $L'=L_x\times\{g\}$. Suppose $\phi'(L')\subset L'$,
which means $\Phi_g(L_x)=L_x$. Any leaf of $\FF'|_{M'\cap(M\times
O)}$ is of the form $\Phi_a(L_x)\times\{a^{-1}ga\}$ for some $a\in
G$. We have
$$
\Phi_{a^{-1}ga}\circ\Phi_a(L_x)=\Phi_{ga}(L_x)=\Phi_a\circ\Phi_g(L_x)=\Phi_a(L_x)\;.
$$
So $\phi'$ preserves $\Phi_a(L_x)\times\{a^{-1}ga\}$.
\end{proof}

According to Proposition~\ref{p:SSS}, if $O_1$ is small enough,
then $\SSS$ is contained in a leaf $M'$ of $\GG$; this property is
assumed from now on. Let $M'_1=M'\cap(M\times O_1)$ and
$\FF'_1=\FF|_{M'_1}$. By Proposition~\ref{p:phi'}, $\phi'$ maps
each leaf of $M'$ to itself, and thus can be restricted to a map
$\phi'_1:M'_1\to M'_1$, which is a foliated map
$(M'_1,\FF'_1)\to(M'_1,\FF'_1)$.

Consider the volume form $\Lambda$ of $G$ as a transverse
invariant measure of $\FF$. By Proposition~\ref{p:M'}-(i),
$\Lambda$ lifts to a transverse invariant measure $\Lambda'_1$ of
$\FF'_1$. Similarly, the Riemannian metric of $M$ lifts to a
Riemannian metric of $M'$, which can be restricted to $M'_1$; the
volume form of this restriction is denoted by $\omega_{M'_1}$.

Even though the foliated manifolds of \cite{HeitschLazarov} are
compact, it is clear that its Lefschetz theorem for foliations
with transverse invariant measures generalizes to the non-compact
case when the transverse invariant measure is compactly supported.

In our case, $M'_1$ may not be compact, but, for every
$f\in\cinfc(O)$ supported in $O_1$,
$\Lambda'_{1,f}=\pr_2^*f\cdot\Lambda'_1$ of $\FF'_1$ is a
compactly supported transverse invariant measure of $\FF'_1$.
Therefore, according to \cite{HeitschLazarov}, the
$\Lambda'_{1,f}$--Lefschetz number $L_{\Lambda'_{1,f}}(\phi'_1)$
of $\phi'_1$ can be defined.

\begin{thm}\label{t:L dis(FF)=L Lambda' N,f(phi' N)}
With the above notation and conditions, we have
\[
\langle L_{\text{\rm
dis}}(\FF),f\rangle=L_{\Lambda'_{1,f}}(\phi'_1)
\]
for every $f\in\cinfc(O)$ supported in $O_1$.
\end{thm}

The proof of Theorem~\ref{t:L dis(FF)=L Lambda' N,f(phi' N)} is
analogous to the proof of Theorem~\ref{t:L dis(FF)=chi Lambda(FF)
cdot delta e}. The holonomy groupoid $\mathfrak{G}'_1$ of $\FF'_1$
can be described like $\mathfrak G$ in Section~\ref{ss:global
action} as a \cinf\ submanifold of $M'_1\times M'_1$ containing
the diagonal. Its penumbras of the diagonal can be also defined
like those of $\mathfrak G$. Its source and target projections are
denoted by $s'_1,r'_1:\mathfrak{G}'_1\to M'_1$. The restriction
$\pr_1\times\pr_1:\mathfrak G'_1\to\mathfrak G$ is a covering map
by Proposition~\ref{p:M'}-(i).

Let $S'_1$ be the \cinf\ vector bundle
$$
s^{\prime*}_1\bigwedge T\FF^{\prime*}_1\otimes r^{\prime*}_1\bigwedge T\FF'_1
$$
over $\mathfrak G'_1$, which can be identified with
$(\pr_1\times\pr_1)^*S$. The space of \cinf\ sections of $S'_1$
supported in penumbras of the diagonal will be denoted by
$C^\infty_\D(S'_1)$. Like in Section~\ref{ss:global action}, there
is a global action of $C^\infty_\D(S'_1)$ in $\Omega(\FF'_1)$.

For each leaf $L'$ of $\FF'_1$, the composite
$\phi^{\prime*}_1\circ\Pi_{L'}$ is a smoothing operator on $L'$, and
let $k'_{\phi,L'}$ denote its smoothing kernel. All of these
smoothing kernels can be combined to define a measurable section
$k_\phi$ of $S'_1$ with \cinf\ restrictions to the fibers of $s'_1$;
$k_\phi$ can be called the {\em leafwise smoothing kernel\/} or {\em
leafwise Schwartz kernel\/} of $\phi^{\prime*}_1\circ\Pi_{\FF'_1}$.
So the restriction of $k_\phi$ to the diagonal $\D_{M'_1}$ is
measurable too. Then $\phi^{\prime*}_1\circ\Pi_{\FF'_1}$ defines an
element of the von~Neumann algebra $W^*(M'_1,\FF'_1,\bigwedge
T\FF_1^{\prime*})$, and we have
\begin{equation}\label{e:L Lambda' 1,f(phi' 1)}
L_{\Lambda'_{1,f}}(\phi'_1)=\Tr_{\Lambda'_1}^s(\phi^{\prime*}_1\circ\Pi_{\FF'_1})
=\int_{M'_1}\Tr^sk_\phi(x,x)\,\omega_{M'_1}\;.
\end{equation}

For any $\psi\in\AAA$ with $\supp\psi\subset[-R,R]$, we have
defined the leafwise Schwartz kernels $k\in\cinf(S)$ and $\tilde
k\in C^\infty_\D(S)$ of $\psi(D_\FF)$ and
$\psi(D_{\widetilde\FF})$ in Section~\ref{s:Gauss-Bonett}.
Similarly, we can define the leafwise Schwartz kernels
$k'_1,k'_\phi\in C^\infty_\D(S'_1)$ of $\psi(D_{\FF'_1})$ and
$\phi^{\prime*}_1\circ\psi(\D_{\FF'_1})$, respectively. It is easy
to see that $k'_1$ can be identified with the lift of $k$ via
$\pr_1\times\pr_1$. Therefore $k'_\phi$ is given by
\begin{equation}\label{e:k' phi}
k'_\phi((x,g),(y,g))=\phi^{\prime*}_1\circ
k'_1(\phi'_1(x,g),(y,g))\equiv\phi_g^*\circ k(\phi_g(x),y)\;.
\end{equation}

Choose a sequence of functions $\psi_m$ like in
Section~\ref{s:Gauss-Bonett}. Let $k$ and $k_{m,t}$ be the
leafwise Schwartz kernels of $\Pi_\FF$ and $\psi_m(t\D_\FF)^2$,
respectively. By \cite[Lemma~1.2]{Roe88II}, $k_{m,t}$ tends to $k$
as $t\to\infty$, and moreover $k_{m,t}$ is uniformly bounded for
large $m$ and $t$. Hence, by~\eqref{e:k' phi}, the leafwise
Schwartz kernel $k'_{\phi,m,t}$ of
$\phi^{\prime*}_1\circ\psi_m(t\D_{\FF'_1})^2$ tends to $k'_\phi$
as $t\to\infty$, and $k_{m,t}$ is uniformly bounded for large $m$
and $t$. Therefore $$
\lim_{t\to\infty}\Tr^s_{\Lambda'_{1,f}}(\phi^{\prime*}_1\circ
\psi_m(t\D_{\FF'_1})^2) =L_{\Lambda'_{1,f}}(\phi'_1) $$ for each
$m$ by~\eqref{e:L Lambda' 1,f(phi' 1)} and the dominated
convergence theorem. Furthermore $$
\Tr^s_{\Lambda'_{1,f}}(\phi^{\prime*}_1\circ
\psi_m(t\D_{\FF'_1})^2) $$ is independent of $t$ (see
\cite[Theorem~5.1]{HeitschLazarov}). Therefore
\begin{equation}\label{e:Tr s=L}
\Tr^s_{\Lambda'_{1,f}}(\phi^{\prime*}_1\circ\psi_m(t\D_{\FF'_1})^2)
=L_{\Lambda'_{1,f}}(\phi'_1)
\end{equation}
for all $m$ and $t$.

Let $Q_{m,t,f}$ be defined like in Section~\ref{s:Gauss-Bonett}.

\begin{lem}\label{l:Tr sQ m,t,f=L Lambda' 0,f(phi' 0)}
We have
$$
\Tr^sQ_{m,t,f}=L_{\Lambda'_{1,f}}(\phi'_1)\;.
$$
\end{lem}

\begin{proof}
By Lemma~\ref{l:psi t}, the Schwartz kernel $q_{m,t,f}$ of
$Q_{m,t,f}$ is given by Corollary~\ref{c:p} when $t$ is small
enough. So, if $(x,x)\in\supp q_{m,t,f}$ for some $x\in M$, we
have $$ q_{m,t,f}(x,x)=\phi_g^*\circ k_{m,t}(\phi_g(x),x)\cdot
f(g)\;, $$ where $g\in O$ is determined by the condition
$x\in\phi_g(B_\FF(x,R_1))$; thus $(x,g)\in\SSS\subset M'_1$.
Therefore, since $\pr_1:\SSS\to M$ is injective
(Proposition~\ref{p:SSS}-(ii)),
\begin{align*}
\Tr^s Q_{m,t,f}& = \int_{\SSS}\Tr^s(\phi_g^*\circ k_{m,t}(\phi_g(x),x))\cdot f(g)\,\omega_{M'_1}(x,g)\\
                         & = \int_{M'_1}\Tr^sk'_{\phi,m,t}((x,g),(x,g))\cdot f(g)\,\omega_{M'_1}(x,g)\\
                         \intertext{by~\eqref{e:k' phi}}
                         & = \Tr^s_{\Lambda'_{1,f}}(\phi^{\prime*}_1\circ\psi_m(t\D_{\FF'_1}))
\end{align*}
for $t$ small enough. Then the result follows by~\eqref{e:Tr s=L}.
\end{proof}

Theorem~\ref{t:L dis(FF)=L Lambda' N,f(phi' N)} follows from
Lemmas~\ref{l:Tr sQ m,t,f=L Lambda' 0,f(phi' 0)} and~\ref{l:lim m
to infty Tr s Q m,t,f}, and Corollary~\ref{c:Tr sQ t,f}.

Now, let us prove Theorem~\ref{t:Lefschetz}. Let $\Fix(\phi')$ and
$\Fix(\phi'_1)$ denote the fixed point sets of $\phi'$ and
$\phi'_1$. Observe that $\Fix(\phi')\subset M'$, and thus
\begin{equation}\label{e:Fix}
\Fix(\phi'_1)=\Fix(\phi')\cap(M\times O_1)\;.
\end{equation}
It is clear that $\pr_2:\Fix(\phi')\to O$ is a proper map because
$M$ is compact and $\Fix(\phi')$ is closed in $M\times O$. Then
$\pr_2:\Fix(\phi'_1)\to O_1$ is proper too by~\eqref{e:Fix}.

A fixed point $(x,g)$ of $\phi'$ is said to be {\em leafwise
simple\/} if $\phi_{g*}-\id : T_x\FF\to T_x\FF$ is an isomorphism.
The set of simple fixed points of $\phi'$ is denoted by
$\Fix_0(\phi')$. Define $\epsilon:\Fix_0(\phi')\to \{\pm 1\}$ by
\[
\epsilon(x,g)=\sign\det (\phi_{g*}-\id : T_x\FF\to T_x\FF)\;.
\]

\begin{lem}\label{l:Fix 0}
$\Fix_0(\phi')$ is a \cinf\ regular submanifold of $M'$ whose
dimension is equal to $\codim\FF$.
\end{lem}

\begin{proof}
Let $\hat\phi:M\times O\to M\times M$ be the \cinf\ map defined by
$\hat\phi(x,g)=(x,\phi_g(x))$, and let $\Delta_M$ denote the
diagonal in $M\times M$. Then
$\Fix(\phi')=\hat\phi^{-1}(\Delta_M)$.

There is some open subset $U\subset M\times O$ such that
$\Fix_0(\phi')=\Fix(\phi')\cap U$. Then the result follows by
showing that the restriction $\hat\phi:U\to M\times M$ is
transverse to $\Delta_M$.

Pick any $(x,g)\in\Fix_0(\phi')$. Let $\Delta_{L_x}$ denote the
diagonal in $L_x\times L_x$. Consider the canonical identity
$T_{(x,x)}(M\times M)\equiv T_xM\oplus T_xM$. The fact that $x$ is
a simple fixed point of $\phi_g$ means that
\begin{equation}\label{e:T x L x oplus T x L x}
T_xL_x\oplus
T_xL_x=\hat\phi_*(T_{(x,g)}(L_x\times\{g\}))+T_{(x,x)}\Delta_{L_x}\;.
\end{equation}
Observe that
\[
\mu_x=\phi_*(T_{(x,g)}(\{x\}\times G))
\]
is complementary of $T_x\FF$, and
\[
\hat\phi_*(T_{(x,g)}(\{x\}\times G))=0_x\oplus\mu_x\;,
\]
where $0_x$ denotes the zero subspace of $T_xM$. So
\begin{align*}
T_xM\oplus T_xM &= (T_xL_x\oplus T_xM)+T_{(x,x)}\Delta_M\\
&=(T_xL_x\oplus T_xL_x)+(0_x\oplus\mu_x)+T_{(x,x)}\Delta_M\\
&=\hat\phi_*(T_{(x,g)}(M\times G))+T_{(x,x)}\Delta_M
\end{align*}
by~\eqref{e:T x L x oplus T x L x}.
\end{proof}

\begin{prop}\label{p:Fix 0}
$\Fix_0(\phi')$ is a $C^\infty$ transversal of $\FF'|_{M'}$.
\end{prop}

\begin{proof}
By Lemma~\ref{l:Fix 0}, it is enough to prove that $\Fix_0(\phi')$
is transverse to $\FF'|_{M'}$, which follows from the following
claim for any point $(x,g)\in\Fix_0(\phi')$.

\begin{claim}\label{cl:T (x,g)(Fix 0(phi)) cap T (x,g) FF'=0}
We have
$$
T_{(x,g)}(\Fix_0(\phi'))\cap T_{(x,g)}\FF'=0\;.
$$
\end{claim}

The proof of Claim~\ref{cl:T (x,g)(Fix 0(phi)) cap T (x,g) FF'=0}
involves another assertion:

\begin{claim}\label{cl:T (x,g)(Fix 0(phi))=ker(phi *-pr 1*}
We have
$$
T_{(x,g)}(\Fix_0(\phi'))=\ker(\phi_*-\pr_{1*}:T_{(x,g)}M'\to T_xM)\;.
$$
\end{claim}

For any $v\in T_{(x,g)}(\Fix_0(\phi'))$, there is a \cinf\ curve
$(x_t,g_t)$ in $\Fix_0(\phi')$, with $-\epsilon<t<\epsilon$ for
some $\epsilon>0$, such that $(x_0,g_0)=(x,g)$ and
$\left.\frac{d}{dt}\,(x_t,g_t)\right|_{t=0}=v$. We have
$\phi(x_t,g_t)=x_t=\pr_1(x_t,g_t)$, yielding
$\phi_*(v)=\pr_{1*}(v)$. So
$$
v\in\ker(\phi_*-\pr_{1*}:T_{(x,g)}M'\to T_xM)\;,
$$
obtaining the inclusion ``$\subset$'' of Claim~\ref{cl:T (x,g)(Fix
0(phi))=ker(phi *-pr 1*}.

Since $\phi_{g*}-\id:T_x\FF\to T_x\FF$ is an isomorphism, so is
$\phi_*-\pr_{1*}:T_{(x,g)}\FF'\to T_x\FF$. Hence
$$
\ker(\phi_*-\pr_{1*}:T_{(x,g)}M'\to T_xM)\cap
T_{(x,g)}\FF'=0\;,
$$
yielding Claims~\ref{cl:T (x,g)(Fix 0(phi)) cap T (x,g)
FF'=0} and~~\ref{cl:T (x,g)(Fix 0(phi))=ker(phi *-pr 1*} because the
inclusion ``$\subset$'' of Claim~\ref{cl:T (x,g)(Fix 0(phi))=ker(phi
*-pr 1*} is already proved.
\end{proof}

\begin{prop}\label{p:submersion}
$\pr_2:\Fix_0(\phi')\to\pr_2(M')$ is a \cinf\ submersion.
\end{prop}

\begin{proof}
Since the leaves of $\FF'$ are contained in the fibers of $\pr_2$,
the tangent map $\pr_{2*}$ induces a homomorphism
$\overline{\pr_{2*}}:T(M\times G)/T\FF'\to TG$. Take any
$(x,g)\in\Fix_0(\phi')$. By Proposition~\ref{p:Fix 0} and
according to the proof of Proposition~\ref{p:GG}, the restrictions
$$
\begin{CD}
T_{(x,g)}\Fix_0(\phi') @>>> T_{(x,g)}M'/T_{(x,g)}\FF' @<<<
\VV_{(x,g)}\;.
\end{CD}
$$
of the quotient map $T(M\times G)\to T(M\times G)/T\FF'$ are
isomorphisms. Moreover $\pr_{2*}$ corresponds to
$\overline{\pr_{2*}}$ by these isomorphisms. So
\begin{align*}
\pr_{2*}(T_{(x,g)}(\Fix_0(\phi')))&=\{(X-\Ad_{g^{-1}}(X))(g))\ |\ X\in\mathfrak{g}\}\\
&=T_g(\pr_2(M'))
\end{align*}
by the proof of Proposition~\ref{p:GG}.
\end{proof}

According to Proposition~\ref{p:Fix 0}, the measure given by
$\Lambda'$ on $\Fix_0(\phi')$ is denoted by
$\Lambda'_{\Fix_0(\phi')}$. The direct image
$\pr_{2*}(\epsilon\cdot\Lambda'_{\Fix_0(\phi')})$ is supported in
$\pr_2(M')\cap O$.

Let $\omega_\Lambda$ be the transverse volume form of $\FF$ defined
by $\Lambda$. Then the transverse volume form of $\FF'|_{M'}$
defined by $\Lambda'$ is $\omega_{\Lambda'}=\pr_1^*\omega_\Lambda$.
The restriction of $\omega_{\Lambda'}$ to the \cinf\ local
transversal $\Fix_0(\phi')$ is a volume form, which can be
identified to the measure $\Lambda'_{\Fix_0(\phi')}$. According to
Proposition~\ref{p:submersion},
$\pr_{2*}(\epsilon\cdot\Lambda'_{\Fix_0(\phi')})$ is given by the
top degree differential form on $\pr_2(M')\cap O$ defined by the
integration along the fibers $$
\int_{\pr_2}\left.\epsilon\cdot\omega_{\Lambda'}\right|_{\Fix_0(\phi')}\;.
$$

By~\eqref{e:Fix}, Theorem~\ref{t:L dis(FF)=L Lambda' N,f(phi' N)}
and the Lefschetz theorem of \cite{HeitschLazarov}, we have
\begin{align*}
\langle L_{\text{\rm dis}}(\FF),f\rangle&=L_{\Lambda'_{1,f}}(\phi'_1)\\
&=\int_{\Fix(\phi'_1)}\epsilon(x,g)\,f(g)\,\Lambda'_{\Fix(\phi')}(x)\\
&=\int_{O_1\cap\pr_2(M')}f(g)\,\pr_{2*}(\epsilon\cdot\Lambda'_{\Fix(\phi')})(g)\\
&=\langle\pr_{2*}(\epsilon\cdot\Lambda'_{\Fix(\phi')}),f\rangle\;,
\end{align*}
completing the proof of Theorem~\ref{t:Lefschetz}.

\section{Examples}

\subsection{Codimension one foliations} Consider the case
when $\FF$ is a codimension one Lie foliation. So we have
${\mathfrak g}=\R$, $G=\R$, and $\FF$ is defined by a closed
nonsingular $1$-form $\omega$. The leaves of $\GG$ in $M\times \R$
are $M'_s=M\times \{s\}$, $s\in\R$. A global $C^\infty$
representation of $\Phi$ is given by the flow $\phi : M\times\R\to
M$ of an arbitrary vector field $X$ on $M$ such that
$\omega(X)=1$. Then
\[
\Fix(\phi')=\{(x,s)\in M\times \R \ |\  \phi_s(x)=x\}\;.
\]
So we have $\Fix(\phi')\cap M'_s\not=\emptyset$ if and only if
either $s=0$ or $s$ is the period of a closed orbit of the flow
$\phi$. In the latter case, we have
\[
\Fix(\phi')\cap M'_s = \bigcup_c\OO_c \times \{s\}\;,
\]
where $c$ runs over the set of all closed orbits of period $s$,
and $\OO_c$ is the corresponding primitive closed orbit:
$$
{\OO}_c=\{\phi_t(x) \in M \ |\  t\in [0,\ell(c)] \}
$$
where $x\in c$ is an arbitrary point, and $\ell(c)$ is the length
of ${\OO}_c$. Assume that all closed orbits of $\phi$ are simple.
Then $\epsilon:\Fix(\phi')\to\{\pm1\}$ is constant on each
${\OO}_c \times \{s\} \subset \Fix(\phi')\cap M'_s$, and its value
on ${\OO}_c \times \{s\}$ will be denoted by $\epsilon_{s}(c)$.

The Lebesgue measure $\Lambda=dt$ on $\R$ can be considered as an
invariant transverse measure of $\FF$. So we have
\[
L_{\text{\rm dis}}(\FF) = \chi_{\Lambda}(\FF)\cdot\delta_0
\]
in some neighborhood of $0$. The restriction of the transverse
volume form $\omega'_\Lambda$ to $\Fix(\phi')\cap M'_s$ coincides
with $\omega_\Lambda$ on each ${\OO}_c$. For any component
${\OO}_c\times \{s\}\subset \Fix(\phi')\cap M'_s$, one can write
$s=k\,\ell(c)$ for some $k\neq0$, and we see that, on
$\R\setminus\{0\}$,
\[
L_{\text{\rm dis}}(\FF) =
\pr_2^*(\left.\epsilon\cdot\Lambda'\right|_{\Fix(\phi')}) =
\sum_c\ell(c)\sum_{k\neq0}\epsilon_{k\,\ell(c)}(c)
\cdot\delta_{k\,\ell(c)}\;,
\]
where $c$ runs over all primitive closed orbits of the flow
$\phi$ \cite{AlvKordy:Betti}.

\subsection{Suspensions} Let $X$ be a connected compact manifold,
$\widetilde{X}$ its universal cover, $G$ a compact Lie group, and
$h:\Gamma=\pi_1(X)\to G$ a homomorphism. Consider the canonical
right action of $\Gamma$ on $\widetilde{X}$, and the diagonal
right action of $\Gamma$ on $\widetilde{M}=\widetilde{X}\times G$:
\[
(x,a)\cdot\gamma=(x\cdot\gamma, h(\gamma^{-1})\cdot a)\;.
\]
Let $M=\widetilde{M}/\Gamma$ (usually denoted by
$\widetilde{X}\times_\Gamma G$). The canonical projection
$\pi:\widetilde{M}\to M$ is a covering map. Let $[x,a]$ be the
element of $M$ represented by each $(x,a)\in \widetilde{M}$. The
foliation $\widetilde{\FF}$ on $\widetilde{M}$ given by the fibers
of the second factor projection $\pr_2: \widetilde{M}\to G$ gives
rise to a foliation $\FF$ on $M$. Let $\Lambda$ be a left invariant
volume form on $G$, which can be considered as an invariant
transverse measure of $\FF$ because its holonomy pseudogroup can be
represented by the pseudogroup generated by the left translations by
elements of $h(\Gamma)$. The corresponding transverse volume form
$\omega_\Lambda$ is defined by the condition
$\pi^*\omega_\Lambda=\pr_2^*\Lambda$ of $\widetilde{\FF}$, whose
restriction to local transversals is another interpretation of
$\Lambda$ as transverse invariant measure of $\FF$. It is easy to
see that
\[
\chi_{\Lambda}(\FF)=\vol (G)\cdot\chi_{\Gamma}(\widetilde{X})\;,
\]
where $\chi_{\Gamma}(\widetilde{X})$ is the $\Gamma$-Euler
characteristic of the covering manifold $\widetilde{X}$ of $X$
defined by Atiyah \cite{Atiyah76}. By Atiyah's $\Gamma$-index
theorem \cite{Atiyah76}, we have $\chi_{\Gamma}(\widetilde{X})=
\chi(X)$, where $\chi(X)$ is the Euler characteristic of $X$.

There is a $C^\infty$ global representation $\phi:M\times G\to M$
of the structural transverse action $\Phi$, defined by
\[
\phi([x,a],g)=[x,ag]\;.
\]
This $\phi$ is a free action. Therefore
\begin{equation}\label{e:suspensions}
L_{\text{\rm dis}}(\FF)=\vol(G)\cdot\chi(X)\cdot\delta_e
\end{equation}
on the whole of $G$. In particular, if $\chi(X)\not=0$, then $\dim
\overline{H}(\FF)=\infty$ for any homomorphism $h:\Gamma\to G$.

We can consider the following concrete example. Let $X$ be a
compact oriented surface of genus $g\ge2$ endowed with a
hyperbolic metric. One can show that there exists an injective
homomorphism $h:\pi_1(X)\to\operatorname{SO}(3,\R)$. One obtains a
Lie $\text{SO}(3,\R)$-foliation $\FF$ whose leaves are dense,
simply connected (diffeomorphic to $\R^2$) and isometric to the
hyperbolic plane. Assuming that $\vol(G)=1$, we get
\[
\beta^0_\Lambda(\FF)=\beta^2_\Lambda(\FF)=0\;,\quad\beta^1_\Lambda(\FF)=2g-2\;.
\]
Since the leaves of $\FF$ are dense, we have
$\overline{H}^0(\R)\cong\overline{H}^2(\R)\cong\R$, and therefore
\[
\Tr^0_{\text{\rm dis}}(\FF)=\Tr^2_{\text{\rm dis}}(\FF)=1\;.
\]
By (\ref{e:suspensions}), we get
\[
L_{\text{\rm dis}}(\FF)=(2-2g) \cdot \delta_e\;,
\]
and
\[
\Tr^1_{\text{\rm dis}}(\FF)=(2g-2) \cdot \delta_e+2\;.
\]

One can also take any homomorphism of $\Gamma$ to the $n$-torus
$\R^n/\Z^n$ to produce a foliation, which has infinite dimensional
reduced cohomology of degree one (see
\cite[Example~2.11]{AlvHector}). In this case, we have
\[
\Tr^i_{\text{\rm dis}}(\FF)\neq\beta^i_{\Lambda}(\FF)\cdot
\delta_e\;,
\]
but $\Tr^i_{\text{\rm dis}}(\FF)-\beta^i_{\Lambda}(\FF)\cdot
\delta_e$ is \cinf.

\subsection{Bundles over homogeneous spaces and the Selberg trace
formula} Let $G$ be a simply connected Lie group, $\Gamma$ a
discrete cocompact subgroup in $G$, and $\alpha$ an injective
homomorphism of $\Gamma$ to the diffeomorphism group $\Diff(X)$ of
some compact connected \cinf\ manifold $X$. Consider a left action
of $\Gamma$ on $\widetilde{M}=G\times X$ given by
\[
\gamma\cdot(a,x)=(\gamma a,\alpha(\gamma)(x))\;.
\]
Let $M=\Gamma\backslash(G\times X)$, and let $[a,x]$ be the
element of $M$ represented by any $(a,x)\in\widetilde{M}$. The
canonical projection $\pi:\widetilde M\to M$ is a covering map.
The first factor projection $G\times X\to G$ defines a fiber
bundle map $M\to \Gamma\backslash G$, whose fibers are the leaves
of a foliation $\FF$. For each $a\in G$, the leaf of $\FF$ through
$\Gamma a$ is
\[
L_{\Gamma a}=\{[a,x]\ |\  x\in X\}\;,
\]
which is diffeomorphic to $X$ because $\alpha$ is injective.
Consider a left-invariant volume form $\Lambda$ on $G$. It induces
a volume form on $\Gamma\backslash G$, denoted by
$\Lambda_{\Gamma\backslash G}$, whose pull-back to $M$ via the map
$M\to \Gamma\backslash G$ defines a transverse volume form
$\omega_\Lambda$ of $\FF$. Since $M\to\Gamma\backslash G$ is a
fiber bundle map with typical fiber $X$, we get
\[
\chi_\Lambda(\FF) = \vol(\Gamma\backslash G)\cdot\chi(X)\;,
\]
where $\chi(X)$ is the Euler characteristic of $X$.

The structural transverse action $\Phi_g$ of an element $g\in G$
is given by the leafwise homotopy class of diffeomorphisms
$\phi_g: M\to M$ of the form
\[
\phi_g([a,x])=[ag,\beta(x)]\;,
\]
where $\beta$ is any diffeomorphism of $X$ homotopic to $\id_X$.

The leaf of the foliation $\GG$ through a point $([a,x],b)\in
M\times G$ is
\[
\{ ([ag, y], g^{-1}bg) \ |\  y\in X,\ g\in G\}\;.
\]
So the leaves of $\GG$ are
\[
M'_b=\{([g,y], g^{-1}bg) \ |\  y\in X,\ g\in G\}\;, \quad b\in
G\;,
\]
with $M'_{b_1}=M'_{b_2}$ when
$b_2\in\operatorname{Ad}(\Gamma)b_1$; thus the leaves of $\GG$ are
parameterized by the $\Gamma$-conjugacy classes in $G$.

Let $\pr_1$ and $\pr_2$ denote the factor projections of $M\times
G$ to $M$ and $G$, respectively. The restriction $\pr_2: M'_b\to
G$ is a bundle map over the orbit
$$
\OO_b=\{g^{-1}bg \ |\  g\in G\}\equiv G_b\backslash G
$$
of the adjoint representation of $G$ on $G$, where
$$
G_b =\{g \in G \ |\  gb=bg\}
$$
is the centralizer of $b$ in $G$.

For each $b\in G$, the restriction $\pr_1: M'_b\to M$ is a
covering map. Indeed, we have $M'_b
\equiv\Gamma_b\backslash(G\times X)$, where
$$
\Gamma_b =\{\gamma\in\Gamma \ |\  \gamma b =b\gamma\}=\Gamma\cap
G_b\;.
$$
The leaves of the foliation $\FF'_b=\pr^*_1\FF$ on $M'_b$ are
described as
\[
L_a=\{([a,y], a^{-1}ba) \ |\  y\in X\}\;, \quad a\in G\;,
\]
with $L_{a_1}=L_{a_2}$ if and only if $\Gamma_ba_1=\Gamma_ba_2$.
Therefore the leaves of $\FF'$ are the fibers of the natural map
$$
M'_b \equiv \Gamma_b\backslash(G\times X) \to\Gamma_b\backslash
G\;, \quad ([a,y], a^{-1}ba) \mapsto \Gamma_b a\;.
$$

Take a $C^\infty$ global representation $\phi:M\times G\to M$ of
$\Phi$ defined by
\[
\phi([a,x], g)=[ag,x]\;.
\]
We have
\[
\Fix(\phi')=\{([a,x], g)\in M\times G \ |\  [ag, x]=[a,x]\}\;.
\]
The identity $[ag, x]=[a,x]$ holds if and only if there exists
$\gamma\in \Gamma$ such that $ag=\gamma a$ and
$\alpha(\gamma)(x)=x$. Hence
\[
\Fix(\phi') =\bigcup_{\gamma\in\Gamma} \{([a,x], a^{-1}\gamma a)\
|\  x\in X,\ \alpha(\gamma)x=x,\ a\in G \}\;.
\]
We see that if $\Fix(\phi')\cap M'_b\neq\emptyset$, then one can
assume that $b=\gamma\in\Gamma$ and $\alpha(\gamma)$ has a fixed
point in $X$. In this case,
\[
\Fix(\phi')\bigcap M'_\gamma=  \{([a,x], a^{-1}\gamma a) \ |\
x\in X,\ \alpha(\gamma)x=x,\ a\in G \}\;.
\]

A point $([a,x],a^{-1}\gamma a)\in \Fix(\phi')\cap M'_\gamma$ is
simple if and only if $x$ is a simple fixed point of
$\alpha(\gamma)$; in this case, we have
\[
\epsilon ([a,x],a^{-1}\gamma a) = \sign \det
(\alpha(\gamma)_*-\id:T_xX\to T_xX)\;,
\]
which is denoted by $\epsilon_{\alpha(\gamma)}(x)$. Assume that, for
any $\gamma\in \Gamma\setminus\{e\}$, all the fixed points of
the diffeomorphism $\alpha(\gamma)$, denoted by $x_1(\gamma),
x_2(\gamma), \dots, x_{d(\gamma)}(\gamma)$, are simple. Then
\[
\Fix(\phi')\bigcap M'_\gamma = \bigcup_{k=1}^{d(\gamma )}
\{([a,x_k(\gamma)], a^{-1}\gamma a) \ |\  a\in G \}\;.
\]
The transverse volume form $\omega'_\Lambda=\pr^*_1\omega_\Lambda$
of $\FF'_\gamma$ is, by definition, the pull-back of
$\Lambda_{\Gamma_\gamma\backslash G}$ via the map $M'_\gamma \to
\Gamma_\gamma\backslash G$. Let $\Sigma$ be a complete set of
representatives of the $\Gamma$-conjugacy classes in $\Gamma$. For $f\in C^\infty_c(G\setminus
\{e\})$, we get
\begin{multline*}
\langle \pr_{2*}(\epsilon\cdot\Lambda'|_{\Fix(\phi')}), f\rangle\\
\begin{aligned}
 & =\int_{\Fix(\phi')} f\circ\pr_2\cdot\epsilon\,\omega'_\Lambda\\
&= \sum_{\gamma\in\Sigma\setminus \{e\}} \sum_{k=1}^{d(\gamma)}
\int_{\Gamma_\gamma\backslash G} f(a^{-1}\gamma a)\cdot
\epsilon_{\alpha(\gamma)}(x_k(\gamma))\,
\Lambda_{\Gamma_\gamma\backslash G}(\Gamma_\gamma a)\;.
\end{aligned}
\end{multline*}
By the classical Lefschetz theorem, we have
\[
\sum_{k=1}^{d(\gamma)} \epsilon_{\alpha(\gamma)}(x_k(\gamma)) =
L(\alpha(\gamma))\;,
\]
where
\[
L(\alpha(\gamma))=\sum_{i=1}^{\dim X}(-1)^i \Tr (\alpha(\gamma)^*
: H^i(X)\to H^i(X))
\]
is the Lefschetz number of the diffeomorphism $\alpha(\gamma)$. It
can be easily seen that $L(\alpha(\gamma))$ depends only on the
conjugacy class of $\gamma$. Take a left invariant Riemannian metric
on $G$ whose volume form is $\Lambda$. Consider the Riemannian
metric on $G_\gamma\backslash G$ so that the canonical projection
$G\to G_\gamma\backslash G$ is a Riemannian submersion, and let
$\Lambda_{G_\gamma\backslash G}$ be the corresponding volume form.
Then
\begin{multline*}
\langle \pr_{2*}(\epsilon\Lambda'|_{\Fix(\phi')}),f\rangle\\
\begin{aligned}
&=\sum_{\gamma\in\Sigma\setminus \{e\}} L(\alpha(\gamma))
\int_{\Gamma_\gamma\backslash G} f(a^{-1}\gamma a)\,
\Lambda_{\Gamma_\gamma\backslash G}(\Gamma_\gamma a)\\
&=\sum_{\gamma\in\Sigma\setminus \{e\}} L(\alpha(\gamma))\,
\vol(\Gamma_\gamma \backslash G_\gamma) \int_{G_\gamma\backslash G}
f(a^{-1}\gamma a)\,\Lambda_{G_\gamma\backslash G}(G_\gamma a)\;.
\end{aligned}
\end{multline*}
Finally, we get the following Selberg type trace formula (cf.
\cite{Selberg}):
\begin{multline*}
\langle L_{\text{\rm dis}}(\FF), f\rangle\\
\begin{aligned}
&= \vol(\Gamma\backslash G)\, \chi(X)\, f(e) \\
& \phantom{=\ }\text{}+ \sum_{\gamma\in\Sigma\setminus \{e\}}
L(\alpha(\gamma))\, \vol(\Gamma_\gamma\backslash G_\gamma)
\int_{G_\gamma \backslash G} f(a^{-1}\gamma a)\, \Lambda_{G_\gamma
\backslash G}(G_\gamma a)\;.
\end{aligned}
\end{multline*}

In the particular case when $G=\R$, $\Gamma =\Z$ and the
homomorphism $\alpha$ is given by a diffeomorphism $F$ of a
compact manifold $X$, the manifold $M$ is the mapping torus of $F$
and the foliation $\FF$ is given by the fibers of the natural map
$M\to S^1$. Then the formula gives
\[
L_{\text{\rm dis}}(\FF) = \chi (X)\cdot\delta_0+
\sum_{k\in\Z\setminus \{0\}}L(F^k)\cdot\delta_k\;.
\]

\subsection{Homogeneous foliations} Let $H$ and $G$ be simply
connected Lie groups, $\Gamma$ a uniform discrete subgroup in $H$,
and $D:H\to G$ a surjective homomorphism so that
$\Gamma_1=D(\Gamma)$ is dense in $G$. Then $M=\Gamma\backslash H$
is a compact manifold,  and let $\FF$ be the foliation on $M$
whose leaves are the projections of the fibers of $D$. If $K=\ker
D$, then the leaves of $\FF$ are the orbits of the right action of
$K$ on $M$ induced by the right action on $H$ defined by right
translations.

This $\FF$ is a Lie $G$-foliation whose structural transverse action
$\Phi$ is given as follows: for each $g\in G$, $\Phi_g$ is
represented by the foliated map $\FF\to \FF$ induced by the right
multiplication by any element of $D^{-1}(g)$.

The leaf of the foliation $\GG$ on $M\times G$ through a point
$(\Gamma h, a)\in M\times G$ is, by definition,
\[
M'_{(\Gamma h, a)}=\{ (\Gamma h_1, g^{-1}ag) \ |\  g\in G,\ h_1
\in D^{-1}(D(\Gamma h)\,g) \}\;.
\]
It is easy to see that there is a bijection between the set of
leaves of $\GG$ and the orbit space $G/\Ad(\Gamma_1)$ of the
adjoint action of $\Gamma_1$ on $G$ so that, for
$\Ad(\Gamma_1)g_0\in G/\Ad(\Gamma_1)$, the corresponding leaf is
described as
\[
M'_{\Ad(\Gamma_1)g_0}=\{(\Gamma h, g)\in M\times G \ |\  D(h)\, g\,
D(h)^{-1}\in \Ad(\Gamma_1)\,g_0\}\;.
\]
The first factor projection $\pr_1: M'_{\Ad(\Gamma_1)g_0}\to M$ is
a covering map; indeed, $M'_{g_0} \equiv \Gamma_{g_0}\backslash
H$, where $\Gamma_{g_0}=\Gamma\cap D^{-1}(\Gamma_{1,g_0})$,
denoting by $\Gamma_{1,g_0}$ the centralizer of $g_0$ in
$\Gamma_1$.

The leaves of $\FF$ can be described as
\[
L_{\Gamma_1 g_1}=\{\Gamma h\in M \ |\  D(h)\in \Gamma_1 g_1\}\;,
\quad g_1\in \Gamma_1\backslash G\;.
\]
By definition, the leaf $L'_{\Gamma_1g_1}=\pr^*_1(L_{\Gamma_1 g_1})$
of the foliation $\FF'=\pr^*_1\FF$ on $M'_{g_0}$ consists of all
$(\Gamma h, g)\in M \times G$ such that $D(h)\, g\, D(h)^{-1}\in
\Ad(\Gamma_1)\,g_0$ and $D(h)\in \Gamma_1 g_1$. So it can be
parameterized by the elements of $(\Gamma_1\backslash G) \times
(G/\Ad(\Gamma_1))$, and it can be described as
\[
L'_{\Gamma_1 g_1} =\{ (\Gamma h, g)\in M \times G : D(h)\in
\Gamma_1 g_1,\ g \in \Ad(g_1)\Ad(\Gamma_1)\,g_0\}\;.
\]

We also see that $\pr_2(M'_{g_0})$ is the orbit $\OO_{g_0}$ of the
adjoint action of $G$ on $G$ through $g_0$. Moreover, $\pr_2 :
M'_{g_0}\to \OO_{g_0}$ is a bundle map, and the fiber of this bundle
over $y\in \OO_{g_0}$ can be identified with $\Gamma_x\backslash
H_x$, where $x\in H$ is any element such that $D(x)=y$.

Denote by $\mathfrak h$, $\mathfrak g$ and $\mathfrak k$ the Lie
algebras of $H$, $G$ and $K$, respectively. We have a short exact
sequence
\[
\begin{CD}
0 @>>> {\mathfrak k} @>>> {\mathfrak h} @>{D_*}>> {\mathfrak g}
@>>> 0\;.
\end{CD}
\]
To construct $C^\infty$ local representations of $\Phi$, we choose a
splitting of this short exact sequence; that is, a linear map $s:
{\mathfrak g} \to {\mathfrak h}$ such that $D_*\circ
s=\id_{\mathfrak g}$. So $s$ is injective and $s({\mathfrak
g})\oplus {\mathfrak k}={\mathfrak h}$. Let $U\subset {\mathfrak g}$
be an open neighborhood of $0$ in ${\mathfrak g}$ such that the
restriction $\exp : U\to \exp(U)\subset G$ of the exponential map to
$U$ is a diffeomorphism. Then, for any $g\in G$, a $C^\infty$ local
representation $\phi: M \times O \to M$ of $\Phi$ is defined on the
open neighborhood $O=g\exp(U)$ of $g$ as
\[
\phi(\Gamma h, g \exp Y)= \Gamma h h_1 \exp s(Y)\;,\quad h\in H\;,
\quad Y\in U\;,
\]
where $h_1\in H$ is any element such that $D(h_1)=g$.

Now fix $g\in G$ and $h_1\in H$ such that $D(h_1)=g$. By
definition,
\[
(\Gamma h, g \exp Y)\in \Fix(\phi') \Leftrightarrow \Gamma h h_1
\exp s(Y) = \Gamma h \Leftrightarrow h h_1 \exp s(Y)\, h^{-1}\in
\Gamma\;.
\]
We have
\[
D(h)\, g\exp Y\, D(h)^{-1}=D(h h_1 \exp s(Y)\, h^{-1})\in \Gamma_1\;,
\]
therefore, we get $\Fix(\phi')\cap M'_{g_0} \not= \emptyset $ iff
$g_0\in\Gamma_1$. In particular, it follows that
\[
\pr_2(\Fix(\phi'))= \bigcup_{\gamma\in\Sigma}\OO_\gamma\;,
\]
where $\Sigma$ is a complete set of representatives of the
$\Gamma_1$-conjugacy classes in $\Gamma_1$. For a fixed class
$\gamma\in \Sigma$, let $[D^{-1}(\gamma)]$ be the
$\Gamma$-conjugacy class of the unique element $\gamma_1\in\Gamma$
such that $D(\gamma_1)=\gamma$. Then we have
\begin{multline*}
\Fix(\phi')\cap M'_{\gamma}\\
=\{ (\Gamma h, g \exp Y)\in
(\Gamma\backslash H) \times G\ |\  h h_1 \exp s(Y)\, h^{-1} \in
[D^{-1}(\gamma)] \}\;.
\end{multline*}

For any $(\Gamma h, g \exp Y)\in \Fix(\phi')$, the left
translation by $h$ determines an isomorphism of the tangent space
$T_{\Gamma h}\FF$ with $\mathfrak k$, and, under this isomorphism,
the induced map $(\phi_{h_1 \exp s(Y)})_* : T_{\Gamma h}\FF \to
T_{\Gamma h}\FF$ corresponds to the restriction $\Ad (h_1 \exp
s(Y))_*\left|_{\mathfrak k}\right. : {\mathfrak k} \to {\mathfrak
k}$ of the differential of the adjoint action of $g \exp Y\in G $
on $G$ to $\mathfrak k$. In particular, $(\Gamma h, g \exp Y)\in
\Fix(\phi)$ is simple if and only if $\Ad (h_1 \exp
s(Y))_*\left|_{\mathfrak k}\right. : {\mathfrak k} \to {\mathfrak
k}$ is an isomorphism. It should be noted that this condition
depends only on $g\exp Y$ and is independent of the choice $h_1$
and $s$.

Assume that $\Ad (h\gamma h^{-1})_*\left|_{\mathfrak k}\right. :
{\mathfrak k} \to {\mathfrak k}$ is an isomorphism for any
$\gamma\in \Gamma$ and $h\in H$. Then the value
\begin{align*}
\epsilon (\Gamma h, g \exp Y) & = \sign \det \left( (\phi_{h_1
\exp s(Y)})_* -\id: T_{\Gamma h}\FF \to T_{\Gamma h}\FF \right)\\
& = \sign \det \left( \Ad (h_1 \exp s(Y))_*\left|_{\mathfrak
k}\right. -\id : {\mathfrak k} \to {\mathfrak k} \right)
\end{align*}
is the same for any $(\Gamma h, g \exp Y)\in \Fix(\phi')\cap
M'_{\gamma}$, and equals
\[
\epsilon(\gamma) = \sign \det \left( \Ad (\gamma)_*
\left|_{\mathfrak k}\right. -\id : {\mathfrak k} \to {\mathfrak k}
\right)\;.
\]

Let $\Lambda$ be a left invariant volume form on $G$, which can be
identified with a transverse volume form of $\FF$. Fix $\gamma\in
\Sigma$. Then the transverse volume form $\Lambda'=\pr_1^*\Lambda$
of $\FF'$ is given by the lift of $\Lambda$ to $M'_{\gamma}$ by the restriction
of the map
\[
(\Gamma h, g)\in (\Gamma\backslash H) \times G \mapsto  D(h) \in
G\;.
\]
to
\[
M'_{\gamma}=\{(\Gamma h, g)\in (\Gamma\backslash H) \times G\ |\
D(h)\, g\, D(h)^{-1}\in \Ad(\Gamma_1)\,\gamma\}\;.
\]
As above, take a left invariant Riemannian metric on $G$ whose
volume form is $\Lambda$. Consider the Riemannian metric on
$G_\gamma\backslash G$ so that the canonical projection $G\to
G_\gamma\backslash G$ is a Riemannian submersion, and let
$\Lambda_{G_\gamma\backslash G}$ be the corresponding volume form.
Restricting the form $\Lambda'$ to $\Fix(\phi')\cap M'_{\gamma}$ and
integrating it along the fibers of $\pr_2$, for any $f\in
C^\infty_c(G)$, we get
\[
\begin{split}
\langle\chi_{\mbox{\scriptsize\rm dis}}(\FF), f\rangle & =\langle
\pr_{2*}(\epsilon\Lambda'), f \rangle\\
& = \sum_{\gamma\in\Sigma}
\epsilon(\gamma) \vol (\Gamma_{\gamma_0}\backslash H_{\gamma_0})
\int_{G_\gamma\backslash G} f(g^{-1}\gamma g)\, \Lambda_{G_\gamma
\backslash G}(G_\gamma g)\;,
\end{split}
\]
where $\gamma_0\in \Gamma $ is the unique element such that
$D(\gamma_0)=\gamma$.

\subsection{Nilpotent homogeneous foliations} Let $G$ be a nontrivial
simply connected nilpotent Lie group and let $\Gamma_1\subset G$ be
a finitely generated dense subgroup. By Malcev's theory
\cite{Mal'cev}, there exists a simply connected nilpotent Lie group
$H$, an embedding $i:\Gamma_1\to H$ and a surjective homomorphism
$D:H\to G$ such that $\Gamma=i(\Gamma_1)$ is discrete and uniform in
$H$, and $D\circ i=\id_{\Gamma_1}$. Consider the corresponding
homogeneous foliation on the closed nilmanifold $M=\Gamma\backslash
H$. As above, $K$ denotes the kernel of $D$, which is a normal
connected Lie subgroup in $H$, and ${\mathfrak k}$ denotes the Lie
algebra of $K$. As shown in \cite[Theorem~2.10]{AlvHector}, there is
a canonical isomorphism $\overline{H}(\FF)\cong H({\mathfrak k})$
({\it c.f.\/} \cite{Nomizu}), and thus $L_{\text{\rm dis}}(\FF)=0$
by Corollary~\ref{c:L=0}. Let us check this triviality in another
way. It can be easily seen that, under this isomorphism, the action
of an element $g\in G$ on $\overline{H}(\FF)$ induced by the
structural action $\Phi$ corresponds to the action $\text{Ad}_*(h)$
on $H({\mathfrak k})$ induced by the adjoint action of any element
$h\in D^{-1}(g)$. So $\Tr^i_{\text{\rm dis}}(\FF)$ is a smooth
function on $G$, whose value at $g\in G$ is the trace of
$\text{Ad}_*(h)$ on $H^i({\mathfrak k})$ with $h\in D^{-1}(g)$.
Since $H$ is nilpotent, $\text{Ad}_*(h)$ has a triangular matrix
representation whose diagonal entries are equal to $1$. So
\[
\Tr^i_{\text{\rm dis}}(\FF)\equiv \dim H^i({\mathfrak k})\;,
\]
yielding
\[
L_{\text{\rm dis}}(\FF)\equiv\sum_i (-1)^i\,\dim H^i({\mathfrak
k}) =\sum_i (-1)^i\,\dim \bigwedge^i{\mathfrak k} = 0\;.
\]
Any local section $g\mapsto h_g$ of $D$ on some open subset
$O\subset G$ induces a \cinf\ local  representation $\phi:M\times
O\to M$ of the structural action $\Phi$, where each $\phi_g$ is
induced by the right multiplication by $h_g$. All the fixed points
of $\phi$ are not simple.

\end{document}